\documentclass{scrartcl}
\usepackage{amsmath}
\usepackage{amsfonts}
\usepackage{amssymb}
\usepackage{dsfont}
\usepackage[T1]{fontenc}
\usepackage[latin1]{inputenc}
\usepackage{epsfig}
\usepackage{graphicx}

\usepackage[boxed, lined]{algorithm2e}
\usepackage{float}
\usepackage{comment}

\newcommand{\bi}{\mathbf{i}}
\newcommand{\ba}{\mathbf{a}}
\newcommand{\bb}{\mathbf{b}}

\newcommand{\D}{{\mathcal{D}}}

\newcommand{\Nu}{{\mathcal{N}}}

\newcommand{\N}{\mathbb{N}}

\newcommand{\R}{\mathbb{R}}

\newcommand{\Rd}{\mathbb{R}^d}

\newcommand{\beq}{\begin{eqnarray*}}

\newcommand{\eeq}{\end{eqnarray*}}

\newcommand{\beqm}{\begin{eqnarray}}

\newcommand{\eeqm}{\end{eqnarray}}

\newtheorem{theorem}{Theorem}

\newtheorem{lemma}{Lemma}
\newtheorem{definition}{Definition}

\newcommand{\EXP}{{\mathbf E}}
\newcommand{\PROB}{{\mathbf P}}
\newcommand{\Var}{{\mathbf{Var}}}

\renewcommand{\bf}{\normalfont \bfseries}
\renewcommand{\it}{\normalfont \itshape}

\allowdisplaybreaks
\begin{document}
\renewcommand{\thefootnote}{\fnsymbol{footnote}}
\newcommand{\F}{{\cal F}}
\newcommand{\Sp}{{\cal S}}
\newcommand{\G}{{\cal G}}
\newcommand{\HH}{{\cal H}}
%\maketitle
%\noindent

\begin{center}

  {\LARGE \bf
    Analysis of the rate of convergence of 
    neural network regression estimates which
    are easy to implement
  }
\footnote{
Running title: {\it Neural network regression estimates}}
\vspace{0.5cm}

Alina Braun$^{1}$, Michael Kohler$^{1}$
and Adam Krzy\.zak$^{2,}$\footnote{Corresponding author. Tel:
  +1-514-848-2424 ext. 3007, Fax:+1-514-848-2830}\\

{\it $^1$
Fachbereich Mathematik, Technische Universit\"at Darmstadt,
Schlossgartenstr. 7, 64289 Darmstadt, Germany,
email: braun@mathematik.tu-darmstadt.de, kohler@mathematik.tu-darmstadt.de}

{\it $^2$ Department of Computer Science and Software Engineering, Concordia University, 1455 De Maisonneuve Blvd. West, Montreal, Quebec, Canada H3G 1M8, email: krzyzak@cs.concordia.ca}

\end{center}
\vspace{0.5cm}

\begin{center}
September 25, 2019
\end{center}
\vspace{0.5cm}

\noindent
    {\bf Abstract}\\
    Recent results in nonparametric regression show that for deep
    learning, i.e., for neural network estimates with many hidden layers,
    we are able to achieve good rates of convergence even in case of
    high-dimensional predictor variables, provided suitable
    assumptions on the structure of the regression function are imposed.
    The estimates are defined by minimizing the empirical
    $L_2$ risk over a class of neural networks. In practice it
    is not clear how this can be done exactly. In this article
    we introduce a new neural network regression estimate where
    most of the weights are chosen regardless of the data
motivated by some recent
approximation results for neural networks, and which is
therefore easy to implement. We show that for this
    estimate we can derive rates of convergence results in case
    the regression function is smooth.
We combine this estimate with the projection pursuit, where we
choose the directions randomly, and we show that 
for sufficiently many repititions
we get a neural network regression
estimate which is easy to implement and which
achieves the one-dimensional rate of convergence 
(up to some logarithmic factor) in case that the regression function
satisfies the assumptions of projection pursuit.

    \vspace*{0.2cm}


\vspace*{0.2cm}


\section{Introduction}
\label{se1}
For many years
neural networks
have been considered as one of the
best approaches in
nonparametric statistics in view of multivariate
statistical applications, in particular in pattern recognition
and in nonparametric regression
(see, e.g.,
the monographs  Hertz, Krogh and Palmer (1991),
Devroye, Gy\"orfi and Lugosi (1996),
Anthony and Bartlett (1999), Gy\"orfi et al. (2002),
Haykin (2008) and Ripley (2008)).
In recent years the focus in applications 
shifted towards so-called deep learning, where
multilayer feedforward neural networks with many hidden
layers are fitted to observed data (see, e.g.,
Schmidhuber (2015) and the literature cited therein).

In this article we study neural network estimates
in the context of nonparametric regression with random
design. Here,
$(X,Y)$ is an $\Rd \times \R$--valued random vector
satisfying $\EXP \{Y^2\}<\infty$, and given a sample
of $(X,Y)$ of size $n$, i.e., given a data set
\begin{equation}
  \label{se1eq1}
\D_n = \left\{
(X_1,Y_1), \dots, (X_n,Y_n)
\right\},
\end{equation}
where
$(X,Y)$, $(X_1,Y_1)$, \dots, $(X_n,Y_n)$ are i.i.d. random variables,
the aim is to construct an estimate
\[
m_n(\cdot)=m_n(\cdot, \D_n):\Rd \rightarrow \R
\]
of the regression function $m:\Rd \rightarrow \R$,
$m(x)=\EXP\{Y|X=x\}$ such that the $L_2$ error
\[
\int |m_n(x)-m(x)|^2 \PROB_X (dx)
\]
is ``small'' (see, e.g., Gy\"orfi et al. (2002)
for a systematic introduction to nonparametric regression and
a motivation for the $L_2$ error).

It is well--known that
one
needs smoothness assumptions on the regression function in
order to derive non--trivial  rate of convergence results
for nonparametric regression estimates
(cf., e.g., Theorem 7.2 and Problem 7.2 in
Devroye, Gy\"orfi and Lugosi (1996) and
Section 3 in Devroye and Wagner (1980)).
To do this we will use the following definition.
\begin{definition}
\label{intde2}
 Let $p=q+s$ for some $q \in \N_0$ and $0< s \leq 1$,
where $\N_0$ is the set of nonnegative integers.
A {\textbf{function}} $f:\R^d \rightarrow \R$ is called
\textbf{$(p,C)$-smooth}, if for every $\alpha=(\alpha_1, \dots, \alpha_d) \in
\N_0^d$
with $\sum_{j=1}^d \alpha_j = q$ the partial derivative
$\frac{
\partial^q f
}{
\partial x_1^{\alpha_1}
\dots
\partial x_d^{\alpha_d}
}$
exists and satisfies
\[
\left|
\frac{
\partial^q f
}{
\partial x_1^{\alpha_1}
\dots
\partial x_d^{\alpha_d}
}
(x)
-
\frac{
\partial^q f
}{
\partial x_1^{\alpha_1}
\dots
\partial x_d^{\alpha_d}
}
(z)
\right|
\leq
C
\cdot
\| x-z \|^s
\]
for all $x,z \in \R^d$, where $\Vert\cdot\Vert$ denotes the Euclidean norm.
\end{definition}
Stone (1982) showed that the optimal minimax rate of convergence in nonparametric
regression for $(p,C)$-smooth functions is $n^{-2p/(2p+d)}$. In case
that $d$ is large compared to $p$ this rate of convergence is
rather slow (so-called curse of dimensionality).
One way to circumvent this curse of dimensionality is to
impose additional constraints on the structure of the regression
function. Stone (1985) assumed that the regression function is
additive, i.e., that $m:\Rd \rightarrow \R$ satisfies
\[
m(x^{(1)},\dots,x^{(d)})
=
m_1(x^{(1)})+ \dots + m_d(x^{(d)})
\quad (x^{(1)},\dots,x^{(d)} \in \R)
\]
for some $(p,C)$--smooth univariate functions $m_1, \dots, m_d:\R
\rightarrow \R$, and showed that in this case suitably defined
spline estimates achieve the corresponding univariate rate of
convergence.
Stone (1994) extended this results to interaction models, where
the regression function is assumed to be a sum of functions
applied to at most $d^* < d$ components of $x$ and showed
in this case that  suitably defined spline estimates
achieve
the $d^*$--dimensional rate of convergence.
Other classes of functions which enable us to achieve a better rate of convergence
results include single index models, where
\[
m(x)=g(a^T x) \quad (x \in \Rd)
\]
for some $a \in \Rd$ and $g:\R \rightarrow \R$
 (cf., e.g.,  H\"ardle and Stoker (1989), H\"ardle, Hall and Ichimura (1993), 
 Yu and Ruppert (2002), Kong and Xia (2007)  and  Lepski and Serdyukova (2014))
and projection pursuit, where
\[
m(x)=\sum_{l=1}^r g_l(a_l^T x) \quad (x \in \Rd)
\]
for some $r \in \N$, $a_l \in \Rd$ and $g_l:\R \rightarrow \R$ $(l=1,
\dots,r)$ (cf., e.g., Friedman and Stuetzle (1981) and Huber (1985)).
In Section 22.3 in Gy\"orfi et al. (2002) it is shown that suitably
defined (nonlinear) least squares estimates achieve in a
$(p,C)$--smooth projection
pursuit model the univariate rate of convergence $n^{-2p/(2p+1)}$ up
to some logarithmic factor.

A generalization of projection pursuit was considered in Horowitz and
Mammen (2007). In this paper the case of a regression function, which satisfies
\[
m(x)=g\left(\sum_{l_1=1}^{L_1}g_{l_1}  \left(\sum_{l_2=1}^{L_2}g_{l_1,l_2}\left( \ldots \sum_{l_r=1}^{L_r}g_{l_1,\ldots, l_r}(x^{l_1,...,l_r}) \right)\right)\right),
\]
where $g, g_{l_1}, \ldots, g_{l_1,\ldots, l_r}$ are $(p,C)$-smooth
univariate functions and $x^{l_1,...,l_r}$ are single components of
$x\in\Rd$ (not necessarily different for two different indices
$(l_1,\ldots,l_r)$), was studied. With the use of a penalized least squares
estimate, the rate $n^{-2p/(2p+1)}$ was proven.

The estimates in Horowitz and Mammen (2007) and the one for projection
pursuit in Section 22.3 in Gy\"orfi et al. (2002) are nonlinear
(penalized) least squares estimates, so in practice it is unclear
how they can be computed exactly.
Friedman and Stuetzle (1981) described easily implementable estimates
for projection pursuit, but in their definition several times
heuristic simplifications are used, and as a consequence it is
unlcear whether for these easily implementable estimates any rate
of convergence result can be shown.

Recently it was shown in several papers that neural networks
can achieve dimensionality reduction in case the regression
function is a composition of (sums of)
functions, where each of the function
is a function of at most $d^*<d$ variables.
The first paper in this respect was
Kohler and Krzy\.zak (2017), where it was shown that
in this case
 suitably defined multilayer
neural networks achieve the rate of convergence  $n^{-2p/(2p+d^*)}$
(up to some logarithmic factor) in case $p \leq 1$.
Bauer and Kohler (2019) showed that this result even holds
for $p>1$ provided the squashing function is suitably
chosen. Schmidt-Hieber (2019) obtained similar results for neural
networks with ReLU activation function, and Kohler and Langer (2019)
showed that the results of Bauer and Kohler (2019) also hold for very simply
constructed fully connected feedforward neural networks.
In Kohler, Krzy\.zak and Langer (2019) it was demonstrated that neural networks
are able to circumvent the curse of dimensionality in case
the regression function has a low local dimensionality.
Results concerning  estimation of regression functions
which are piecewise polynomials
with partitions with rather general smooth boundaries
by neural networks
 have been
derived in Imaizumi and Fukamizu (2019).

In all articles above the neural network regression estimate is defined
as a nonlinear least squares estimate.
For instance, 
an estimate is defined as the function $m_n \in \F_n$ which minimizes
the empirical $L_2$ risk
\begin{equation}
  \label{se1eq2}
  \frac{1}{n}
  \sum_{i=1}^n
  |Y_i-m_n(X_i)|^2
  \end{equation}
over a nonlinear class $\F_n$ of neural
networks. In practice, it is usually not possible to find the
global minimum of the empirical $L_2$ risk over a class of neural networks and 
one usually tries to find a local minimum using, e.g., the
gradient descent algorithm (so--called backpropagation).

There exist quite a few papers which try to show that neural network
estimates learned by backpropagation have nice theoretical properties.
The most popular approach in this context is the so--called
landscape approach. 
Choromanska et al. (2015)
used random matrix theory to
derive a heuristic argument showing
that the  risk  of  most  of  the  local  minima  of the
empirical $L_2$ risk is  not  much
larger  than  the  risk  of  the  global  minimum. For neural
networks with special activation function it was possible to 
validate this claim, see, e.g., Arora et al. (2018), 
Kawaguchi (2016),
and Du and Lee (2018), which have 
analyzed gradient descent for neural networks with linear or quadratic
activation function. 
But for such neural networks there do not exist good approximation results, 
consequently, one cannot derive from these results 
rates of convergence comparable to the ones above for the
least squares neural network regression estimates.

Du et al. (2018) 
analyzed gradient descent applied to neural networks
with one hidden layer in case of an input with a Gaussian distribution. 
They used the expected gradient instead of the gradient 
in their gradient descent routine, 
and therefore, 
their result cannot be used to derive a rate of convergence result
similar to the results for the least squares neural network estimates
cited above
for an estimate learned by the gradient descent.
Liang et al. (2018)
applied gradient descent to a modified loss function in classification, 
where it is assumed that the data can be interpolated by a neural network.
Here, the last assumption is not satisfied in nonparametric regression
and it is unclear whether the main idea (of simplifying the estimation by a modification of the loss function) 
can also be used in regression setting.
Recently it was shown in several papers, see, e.g.,
Allen-Zhu, Li and Song (2019), Kawaguchi and Huang (2019)
and the literature cited therein, that gradient descent leads to 
a small empirical $L_2$ risk in over-parametrized neural networks. 
Here, it is unclear what the $L_2$ risk of the estimate is (and a bound on this term is necessary in order to derive results like the ones cited above
for the least squares neural network regression estimates).
In particular, due to the fact that the networks are over-parametrized,
a bound on the empirical $L_2$ risk might be not useful for bounding the $L_2$ risk.
And the bound on the $L_2$ risk presented in Kawaguchi and Huang
(2019) requires that the weights in the network be small, and it is
not clear whether this will be satisfied in an over-parametrized neural
network
learned by the gradient descent.

So although the above theoretical results
for the least squares neural network estimates
are quite impressive, there is a big gap between the
estimates for which one has proven the above mentioned
nice rate of convergence results
and the estimates which can be used in
practice. And until now, the results derived in the literature
for backpropagation are unfortunately not strong enough to narrow
this gap.

In this paper we are interested in the following question:
If we define a neural network regression estimate theoretically
exactly as it is implemented in practice, what 
rate of convergence result can we show for this estimate?
This question was already considered in Braun, Kohler and Walk (2019).
There neural network regression estimates
with one hidden layer have been considered, where the weights
were chosen by minimizing a regularized empirical $L_2$ risk
via backpropagation with starting values chosen repeatedly
randomly from a special structure adopted to projection
pursuit. It was shown in a $(p,C)$--smooth projection pursuit
model, i.e., in a projection pursuit model with $(p,C)$--smooth functions, 
that this easily implementable
estimate achieves (up to a logarithmic factor) the rate of convergence
$n^{-2p/(2p+1)}$, provided  $p \leq 1$.
%Here the estimate
%can be computed in time $O(n^2)$ if the number
%$r$ of terms in the projection pursuit model satisfies
%$r \cdot d \leq 2p+1$ (so for $d>1$ we need to have $r=1$).

In the sequel we use a different (but related) approach in order
to derive rate of convergence results for easily implementable
neural network estimates. We use neural
networks with several hidden layers where
we choose most of the inner weights of the network in a
data-independent way, and where we choose the weights of the
output level via regularized least squares estimates.
Here the choice of the inner weights is motivated by recent
approximation results derived for deep neural networks,
and the use of the regularized least squares criterion
for the choice of the weights of the output layer
leads to estimates
which are easy to implement (because they can be computed
by solving a linear equation system).

Our first main result is that in this way we can define
neural network regression estimates which are easy to implement
and which achieve the same rate of convergence results as linear
regression estimates (e.g., kernel or spline estimates), i.e.,
they achieve (up to some logarithmic factor) the optimal
minimax rate of convergence $n^{-2p/(2p+d)}$ in case of a
$(p,C)$--smooth regression function, for any $p >0$.
For our second main result we define in a projection pursuit
model a neural network regression estimate, where
we choose the directions of this model several times randomly
and define the inner weights independent of the data using these random
directions, and where the weights of the output layer are 
computed by using a regularized least squares criterion.
For this estimate we show
that for sufficiently many repetitions (of the choices of the
random directions)
we get an estimate which achieves the one-dimensional rate of convergence 
(up to some logarithmic factor)
in case that the regression function
satisfies the assumptions of the projection pursuit.
To our knowledge this result is the first result
in the literature which shows that there exist
estimates 
which can be easily implemented and which achieve
(up to a logarithmic factor)
the rate of convergence
$n^{-2p/(2p+1)}$ in a $(p,C)$--smooth projection pursuit model
for arbitary $p>0$.

Throughout the paper, the following notation is used:
The sets of natural numbers, natural numbers including $0$,
and real numbers
are denoted by $\N$, $\N_0$ and $\R$, respectively. For $z \in \R$, we denote
the smallest integer greater than or equal to $z$ by
$\lceil z \rceil$.
Furthermore we set $z_+=\max\{z,0\}$.
The  Euclidean norm of $x \in \Rd$
is denoted by $\|x\|$ and $\|x\|_\infty$
denotes its supremum norm.
For $f:\R^d \rightarrow \R$ 
\[
\|f\|_\infty = \sup_{x \in \R^d} |f(x)|
\]
is its supremum norm. %, and the supremum norm of $f$
%on a set $A \subseteq \R^d$ is denoted by
%\[
%\|f\|_{\infty,A} = \sup_{x \in A} |f(x)|.
%\]
Let $\F$ be a set of functions $f:\Rd \rightarrow \R$,
let $x_1, \dots, x_n \in \Rd$ and set $x_1^n=(x_1,\dots,x_n)$.
A finite collection $f_1, \dots, f_N:\Rd \rightarrow \R$
  is called an $\varepsilon$-- cover of $\F$ on $x_1^n$
  if for any $f \in \F$ there exists  $i \in \{1, \dots, N\}$
  such that
  \[
\frac{1}{n} \sum_{k=1}^n |f(x_k)-f_i(x_k)| < \varepsilon.
  \]
  The $\varepsilon$--covering number of $\F$ on $x_1^n$
  is the  size $N$ of the smallest $\varepsilon$--cover
  of $\F$ on $x_1^n$ and is denoted by $\Nu_1(\varepsilon,\F,x_1^n)$.

The outline of this paper is as follows: In Section \ref{se2} the
newly proposed neural network regression estimates for $(p,C)$--smooth
regression functions are defined and a
result for the rate of convergence of these estimates is
presented. In Section \ref{se3} we describe how these
estimates can be combined with projection pursuit, and present
a rate of convergence result where the easily computable estimate
achieves (up to some logarithmic factor) the
optimal one-dimensional rate of convergence if the regression
function satisfies the assumptions of projection pursuit.
The finite sample size performance of our newly proposed estimates
on simulated data
is illustrated in Section \ref{se4}. The proofs
are given in Section \ref{se5}.

\section{A linear (regularized) least squares regression estimate}
\label{se2}
The starting point in defining a neural network is the
choice of an activation function $\sigma: \mathbb{R} \to \mathbb{R}$.
Here, we use in the sequel
so--called squashing functions, which are nondecreasing
and satisfy $\lim_{x \rightarrow - \infty} \sigma(x)=0$
and
$\lim_{x \rightarrow  \infty} \sigma(x)=1$.
An example of a squashing function is
the so-called sigmoidal or logistic squasher
\begin{equation}
  \label{se2eq4}
\sigma(x)=\frac{1}{1+\exp(-x)} \quad (x \in \R).
\end{equation}
The network architecture $(L, \textbf{k})$ depends on a positive
integer $L$ called the \textit{number of hidden layers} and a
\textit{width vector} $\textbf{k} = (k_1, \dots, k_{L}) \in
\mathbb{N}^{L}$ that describes the number of neurons in the first,
second, $\dots$, $L$-th hidden layer. A multilayer feedforward neural
network with
architecture $(L, \textbf{k})$ and sigmoidal function $\sigma$ is a real-valued function 
%defined on $\mathbb{R}^d$ of the form
$f: \R^d \rightarrow \R$ defined by
\begin{equation}\label{se2eq1}
f(x) = \sum_{i=1}^{k_L} c_i^{(L)} \cdot f_i^{(L)}(x) + c_0^{(L)}
\end{equation}
for some $c_0^{(L)}, \dots, c_{k_L}^{(L)} \in \mathbb{R}$ and for $f_i^{(L)}$'s recursively defined by
\begin{equation}
  \label{se2eq2}
f_i^{(r)}(x) = \sigma\left(\sum_{j=1}^{k_{r-1}} c_{i,j}^{(r-1)} \cdot f_j^{(r-1)}(x) + c_{i,0}^{(r-1)} \right)
\end{equation}
for some $c_{i,0}^{(r-1)}, \dots, c_{i, k_{r-1}}^{(r-1)} \in \mathbb{R}$
$(r=2, \dots, L)$
and
\begin{equation}
  \label{se2eq3}
f_i^{(1)}(x) = \sigma \left(\sum_{j=1}^d c_{i,j}^{(0)} \cdot x^{(j)} + c_{i,0}^{(0)} \right)
\end{equation}
for some $c_{i,0}^{(0)}, \dots, c_{i,d}^{(0)} \in \mathbb{R}$.

In the sequel we want to use the data (\ref{se1eq1}) in order
to choose the weights of the neural network such that the resulting
function defined by (\ref{se2eq1})--(\ref{se2eq3}) is a good
estimate of the regression function. To do this, we will fix the
network architecture of the neural network and all weights except
the weights in the output layer and will use the data
 (\ref{se1eq1})
together with the principle of (regularized) least squares in order
to estimate the weights in the output layer.

\subsection{Definition of the network architecture}
\label{se2sub1}
Let $a>0$ be fixed.
The choice of the
network architecture and of the values
values of most of the weights of our neural network
is motivated by the following
approximation result of a $(p,C)$-smooth function for $x \in [-a,a]^d$ by
a local convex
combination of Taylor polynomials:
For
$M \in \N$
and
$\bi= (i^{(1)}, \dots, i^{(d)}) \in \{0, \dots, M\}^d$ set
\[
x_{\bi}
=\left(-a + i^{(1)} \cdot \frac{2a}{M}, \dots, -a + i^{(d)} \cdot \frac{2a}{M}
\right)
\]
and let
\[
\{ \bi_1, \dots, \bi_{(M+1)^d}\}
=
\{0, \dots, M\}^d.
\]
For $k \in \{1, \dots, (M+1)^d\}$ let
\begin{equation*}
p_{\bi_k}(x) = \sum_{\substack{j_1, \dots, j_d \in \{0, \dots, q\}\\j_1 + \dots+j_d \leq q}} \frac{1}{j_1! \cdots j_d!} \cdot \frac{\partial^{j_1+\dots+j_d}f}{\partial^{j_1}x^{(1)} \cdots \partial^{j_d}x^{(d)}}(x_{\bi_k}) \cdot (x^{(1)}-x_{\bi_k}^{(1)})^{j_1} \cdots (x^{(d)}-x_{\bi_k}^{(d)})^{j_d}
\end{equation*}
be the the Taylor polynomial of $f$ with order $q$ around $x_{\bi_k}$ and set
\begin{equation}\label{se2eq5}
P(x) = \sum_{k=1}^{(M+1)^d} p_{\bi_k}(x) \prod_{j=1}^d \left( 1-\frac{M}{2a}
\cdot |x^{(j)} - x_{\bi_k}^{(j)}| \right)_+,
\end{equation}
where $z_+=\max\{z,0\}$ $(z \in \R)$.
Since $P(x)$ is a local convex combination of Taylor polynomials of
$m$,
it is possible to show that for a $(p,C)$--smooth function $m$
we have
\begin{equation}
  \label{se2eq5b}
\sup_{x \in [-a,a]^d}
|m(x)-P(x)|
\leq
c_1
\cdot \frac{1}{M^p}
\end{equation}
(cf.,  Lemma 5 in Schmidt--Hieber (2019)).

$P(x)$
can be written in the form
\[
\sum_{k=1}^{(M+1)^d}
\sum_{\substack{j_1, \dots, j_d \in \{0, \dots, q\}\\j_1 + \dots+j_d \leq q}}
a_{\bi_k,j_1,\dots,j_d}
\cdot (x^{(1)}-x_{\bi_k}^{(1)})^{j_1} \cdots (x^{(d)}-x_{\bi_k}^{(d)})^{j_d}
\prod_{j=1}^d \left(1-\frac{M}{2a} \cdot |x^{(j)} - x_{\bi_k}^{(j)}| \right)_+
\]
with appropriately chosen $a_{\bi_k ,j_1,\dots,j_d} \in \R$.
Our main trick in the sequel is to define
appropriate neural networks $f_{net,j_1,\dots,j_d,\bi_k}$ which approximate
the functions
\[
x \mapsto
(x^{(1)}-x_{\bi_k}^{(1)})^{j_1} \cdots (x^{(d)}-x_{\bi_k}^{(d)})^{j_d}
\prod_{j=1}^d(1-\frac{M}{2a} \cdot |x^{(j)} - x_{\bi_k}^{(j)}|)_+,
\]
and to choose the network architecture such that neural networks
of the form
\[
\sum_{k=1}^{(M+1)^d}
\sum_{\substack{j_1, \dots, j_d \in \{0, \dots, q\}\\j_1 + \dots+j_d \leq q}}
a_{\bi_k,j_1,\dots,j_d}
\cdot
f_{net,j_1,\dots,j_d,\bi_k}(x)
\quad
(a_{\bi_k ,j_1,\dots,j_d} \in \R)
\]
are contained in it. To do this, we let
$\sigma(x)=1/(1+\exp(-x))$ $(x \in \R)$ be the logistic squasher,
choose $R \geq 1$
and define the following neural networks:
The neural network
\begin{equation}
\label{se2eq6}
f_{id}(x)=
4 R \cdot
\sigma \left(
\frac{x}{R}
\right)
-
2R
\end{equation}{
which approximates the function $f(x)=x$ (cf., Lemma \ref{le1} below), the neural network
\begin{eqnarray}
\label{se2eq7}
f_{mult}(x,y)
&=&
\frac{R^2}{4}
\cdot
\frac{(1+e^{-1})^3}{e^{-2}-e^{-1}}
\cdot
\Bigg(
\sigma \left(
\frac{2 (x+y)}{R}+1
\right)
-
2 \cdot \sigma \left(
\frac{x+y}{R} +1
\right)
\nonumber
\\
&&
\hspace*{3cm}
-
\sigma \left(
\frac{2(x-y)}{R}+1
\right)
+
2 \cdot
\sigma \left(
\frac{x-y}{R}+1
\right)
\Bigg),
\end{eqnarray}
which approximates the function $f(x,y)=x \cdot y$
(cf., Lemma \ref{le2} below), the neural network
\begin{equation}
\label{se2eqReLu}
f_{ReLu}(x)
=
f_{mult}(f_{id}(x),\sigma(R \cdot x))
\end{equation}
which approximates $f(x)=x_+$
(cf., Lemma \ref{le3} below), and the neural network
\[
f_{hat,y}(x)
=
f_{ReLU} \left(
\frac{M}{2a} \cdot (x-y)+1
\right)
-
2
\cdot
f_{ReLU} \left(
\frac{M}{2a} \cdot (x-y)
\right)
+
f_{ReLU} \left(
\frac{M}{2a} \cdot (x-y) -1
\right)
\]
which approximates for fixed $y \in \R$ the function
$f(x)=(1-(M/(2a)) \cdot|x-y|)_+$
(cf., Lemma \ref{le4} below).

With these networks we can now define
$f_{net,j_1,\dots,j_d,\bi_k}$
recursively as follows:
We choose $N \geq q$, set $s=\lceil \log_2(N+d) \rceil$ and
define for $j_1, \dots, j_d \in \{0,1,\dots,N\}$
and
$k \in \{1, \dots, (M+1)^d\}$
\[
f_{net,j_1,\dots,j_d,\bi_k}(x)=f_1^{(0)}(x),
\]
where
\[
f_k^{(l)}(x)
=
f_{mult}
\left(
f_{2k-1}^{(l+1)}(x),f_{2k}^{(l+1)}(x)
\right)
\]
for $k \in \{1, 2, \dots, 2^l\}$ and $l \in \{0,\dots,s-1\}$,
and
\[
f_k^{(s)}(x)=f_{id}(f_{id}(x^{(l)}-x_{\bi_k}^{(l)}))
\]
for $j_1 + j_2 + \dots + j_{l-1}+1 \leq k \leq j_1 + j_2 + \dots +
j_l$
and $l=1,\dots,d$,
\[
f_{j_1+j_2+\dots+j_d+k}^{(s)}(x) = f_{hat,x_{\bi_k}^{(k)}}(x^{(k)})
\]
for $k=1, \dots, d$, and
\[
f^{(s)}_k(x)=1
\]
for $k=j_1+j_2+\dots+j_d+d+1, j_1+j_2+\dots+j_d+d+2, \dots,2^s$.
It is easy to see that
$f_{net,j_1,\dots,j_d,\bi_k}$ is a neural network with $s+2$ hidden
layers
and at most
\[
6 \cdot 2^s, \,
12 \cdot 2^s, \,
2 \cdot 2^s,2^s,\,
\dots,
8,4
\]
neurons in the layers $1,2,\dots,s+2$, resp.
Consequently, this
 network
is contained in the class of all fully connected neural networks
with $s+2$ hidden layers and $24 \cdot (N+d)$ neurons
in each hidden layer.
Furthermore it is easy to see that all weights are bounded in absolute
value by $c_2 \cdot \max\{ 1, M/a, R^2 \}$.

\subsection{Definition of the output weights}
\label{se2sub2}
We define our neural network estimate $\tilde{m}_n(x)$ by
\[
\tilde{m}_n(x)
=
\sum_{k=1}^{(M+1)^d}
\sum_{\substack{j_1, \dots, j_d \in \{0, \dots, N\}\\j_1 + \dots+j_d \leq N}}
a_{\bi_k,j_1,\dots,j_d}
\cdot
f_{net,j_1,\dots,j_d,\bi_k}(x),
\]
where the coefficients $a_{\bi_k ,j_1,\dots,j_d}$ are chosen by minimizing
\begin{equation}
\label{se2sub2eq1}
\frac{1}{n} \sum_{i=1}^n |Y_i-\tilde{m}_n(X_i)|^2
+
\frac{c_3}{n} \cdot
\sum_{k=1}^{(M+1)^d}
\sum_{\substack{j_1, \dots, j_d \in \{0, \dots, N\}\\j_1 + \dots+j_d \leq N}}
a_{\bi_k,j_1,\dots,j_d}^2
\end{equation}
for some constant $c_3>0$.
This regularized linear least squares estimate can be computed by
solving a linear equation system. To see this, set
\[
J= (M+1)^d \cdot \left({N+d \atop d}  \right),
\]
let
\[
\{ B_{j} \, : \, j=1,\dots, J\}
=
\left\{
f_{net,j_1,\dots,j_d,\bi_k}(x)
\, : \,
1 \leq k \leq (1+M)^d
\mbox{ and }
 0 \leq j_1 + \dots+j_d \leq N
\right\}
\]
and set
\[
{\mathbf B} = (B_j(X_i))_{1 \leq i \leq n , 1 \leq j \leq J}
\quad
\mbox{and}
\quad
{\mathbf Y}=(Y_i)_{i=1,\dots,n}.
\]
It is easy to see
(cf., Supplement for a corresponding proof)
that the vector of coefficients of our estimate 
is the unique solution of the linear equation system
\begin{equation}
\label{se2eq9}
\left(
\frac{1}{n}
 {\mathbf B}^T  {\mathbf B} +
\frac{c_3}{n} \cdot {\mathbf 1}
\right)
{\mathbf a}
=
\frac{1}{n}
{\mathbf B}^T {\mathbf Y}.
\end{equation}
The value of (\ref{se2sub2eq1})
will be
also
less than or equal to
the value which we get for coefficients equal to zero,
hence we have
\[
\frac{1}{n} ({\mathbf Y}- {\mathbf B a})^T ({\mathbf Y} - {\mathbf B
  a})
+
\frac{c_3}{n} \cdot{\mathbf a}^T {\mathbf a}
\leq
\frac{1}{n} \sum_{i=1}^n Y_i^2,
\]
which will allow us to derive a bound on the maximal absolute value
of our coefficients.

\subsection{Rate of convergence}
\label{se2sub3}

\begin{theorem}
  \label{th1}
  		Assume that the distribution of $(X,Y)$ satisfies
		\begin{align}\label{subgaus}
		\mathbf E \left( e^{c_{4}\cdot |Y |^2}\right) <\infty
		\end{align}
		for some constant $c_{4}>0$ and that the distribution of $X$ has bounded
		support $supp(X)$, and let $m(x)=\EXP\{Y|X=x\}$ be the corresponding
		regression function.
                Assume that $m$ is $(p,C)$--smooth, where
                $p=q+s$ for some $q \in \N_0$
and $s \in (0,1]$.
Define the estimate $\tilde{m}_n$ as in Subsection \ref{se2sub2},
where $\sigma$ is the logistic squasher and where
$N \geq q$,
$M=M_n=\lceil c_5 \cdot n^{1/(2p+d)} \rceil$,
$R=R_n=n^{d+4}$
and
$a=a_n=(\log n)^{1/(6(N+d)}$.
Set
$\beta_n = c_{6} \cdot \log(n)$ for some suitably large constant $c_{6}>0$
and define $m_n$ by
\[
m_n(x)=T_{\beta_n} \tilde{m}_n(x)
\]
where $T_\beta z = \max\{ \min\{z,\beta\},-\beta\}$
for $z \in \R$ and $\beta>0$.
Then
 $m_n$ satisfies for
$n$ suffciently large
\[
  \mathbf E \int |m_n(x) - m(x)|^2 \PROB_X (dx) \leq
c_7 \cdot (\log n)^3 \cdot n^{-\frac{2p}{2p+d}},
\]
where $c_7 >0$ does not depend on $n$.
  \end{theorem}

\noindent
    {\bf Remark 1.}
    It follows from the proof of Theorem \ref{th1}
    that the result also holds for more general squashing
    functions than the logistic squasher. More precisely, in case that the
    definitions of $f_{id}$, $f_{mult}$ and $f_{ReLU}$
    are modified as in Lemma \ref{le1}, Lemma \ref{le2} and Lemma \ref{le3} below,
    it suffices to assume that
    $\sigma$ is Lipschitz continuous and $2$--admissible according to Definition \ref{se5de1}
    below.

%\section{An easily computable neural network regression estimate}
\section{Extension to projection pursuit}
\label{se3}

In this section we assume that
the regression function satisfies
\[
m(x)=\sum_{l=1}^r g_l \left(
a_{(l-1) \cdot d+1} \cdot x^{(1)}
+
\dots
+
a_{l\cdot d} \cdot x^{(d)}
\right)\quad (x^{(1)}, \dots, x^{(d)} \in \R)
\]
for some $r \in \N$, 
some $(p,C)$--smooth functions
$g_l:\R \rightarrow \R$ $(l = 1, \dots, r)$
and some $\ba_l =(a_{(l-1) \cdot d +1}, \dots, a_{l \cdot d})^T \in \R^d$
with $\|\ba_l\|=1$
$(l=1, \dots, r)$. 
Our goal is to construct a neural network
regression estimate of $m$ which achieves the univariate rate
of convergence. 

\subsection{Definition of the network architecture}
\label{se3sub1}
Let $A>0$ be fixed.
The choice of the
network architecture and of the values 
of most of the weigths of our neural network
is motivated by the following
approximation result
for $x\in [-A,A]^d$: 
For $M \in \N$ and $i \in \{0,\dots,M\}$
set
\[
u_i=-\sqrt{d} \cdot A + i \cdot
\frac{2 \cdot \sqrt{d} \cdot A}{M}
%\quad (k=0, \dots, K).
\]
and let
\[
\{i_1, \dots, i_{M+1}\} = \{0, \dots, M\}.
\]
We will see in Section \ref{se5} below that we can approximate
a $(p,C)$--smooth projection pursuit model
\[
m(x) = \sum_{l=1}^r g_l(\ba_l^T x)
\]
by choosing $\bb_l$ close to $\ba_l$ and by
choosing an appropriate
sum of local convex combinations
of polynomials of the form
\[
\sum_{l=1}^r
\sum_{k=1}^{M+1}
\sum_{j_1, \dots, j_d \in \{0,\dots,q\}, \atop j_1+\dots+j_d \leq q}
a_{i_k,j_1,\dots,j_d,\bb_l}
\cdot
(x^{(1)})^{j_1}
\cdot
\dots
\cdot
(x^{(d)})^{j_d}
\cdot
\left(
1
-
\frac{M}{2 \cdot \sqrt{d} \cdot A}
\cdot
|\bb_l^T x - u_{i_k}|
\right)_+.
\]
Our main trick in the sequel is to define
appropriate neural networks $f_{net,j_1,\dots,j_d,i_k,\bb_l}$ which approximate
the functions
\[
x \mapsto
(x^{(1)})^{j_1} \cdots (x^{(d)})^{j_d}
\cdot
\left(
1
-
\frac{M}{2 \cdot \sqrt{d} \cdot A}
\cdot
|\bb_l^T x - u_{i_k}|
\right)_+
\]
and to choose the network architecture such that neural networks
of the form
\[
\sum_{l=1}^r
\sum_{k=1}^{M+1}
\sum_{\substack{j_1, \dots, j_d \in \{0, \dots, q\}\\j_1 + \dots+j_d \leq q}}
a_{i_k,j_1,\dots,j_d,\bb_l}
\cdot
f_{net,j_1,\dots,j_d,i_k,\bb_l}(x)
\quad
(a_{i_k ,j_1,\dots,j_d,\bb_l} \in \R)
\]
are contained in it. To do this, we let
$\sigma(x)=1/(1+\exp(-x))$ $(x \in \R)$ be the logistic squasher,
choose $R \geq 1$
and define the following neural networks:
The neural network
$
f_{id}(x)
$
as in $(\ref{se2eq6})$
	which approximates the function $f(x)=x$ (cf., Lemma \ref{le1} below), the neural network
	$
	f_{mult}(x,y)
	$
	as in $(\ref{se2eq7})$
	which approximates the function $f(x,y)=x \cdot y$
	(cf., Lemma \ref{le2} below), the neural network
	$
	f_{ReLu}(x)
	$
	as in $(\ref{se2eqReLu})$
	which approximates $f(x)=x_+$
	(cf., Lemma \ref{le3} below), and the neural network
	\begin{eqnarray*}
		\bar{f}_{hat,y}(x)
		&=&
		f_{ReLU} \left(
		\frac{M}{2 \cdot \sqrt{d} \cdot A}
		\cdot (x-y)+1
		\right)
		-
		2
		\cdot
		f_{ReLU} \left(
		\frac{M}{2 \cdot \sqrt{d} \cdot A}
		\cdot (x-y)
		\right)
		\\
		&&
		+
		f_{ReLU} \left(
		\frac{M}{2 \cdot \sqrt{d} \cdot A}
		\cdot (x-y) -1
		\right)
	\end{eqnarray*}
	which approximates for fixed $y \in \R$ the function
	$f(x)=(1 - \frac{M}{2 \cdot \sqrt{d} \cdot A}
	\cdot|x-y|)_+$
	(cf., Lemma \ref{le4} below).
	
	With these networks we can now define
	$f_{net,j_1,\dots,j_d,i_k,\bb_l}$
	recursively as follows:
	We choose $N \geq q$, set $s=\lceil \log_2(N+1) \rceil$ and
	define for $l \in \{1, \dots,r\}$,
	$j_1, \dots, j_d \in \{0,1,\dots,N\}$
	and
	$k \in \{1, \dots, M+1\}$
	\[
	f_{net,j_1,\dots,j_d,i_k,\bb_l}(x)=f_1^{(0)}(x),
	\]
	where
	\[
	f_k^{(t)}(x)
	=
	f_{mult}
	\left(
	f_{2k-1}^{(t+1)}(x),f_{2k}^{(t+1)}(x)
	\right)
	\]
	for $k \in \{1, 2, \dots, 2^t\}$ and $t \in \{0,\dots,s-1\}$,
	and
	\[
	f_k^{(s)}(x)=f_{id}(f_{id}(x^{(t)}))
	\]
	for $j_1 + j_2 + \dots + j_{t-1}+1 \leq k \leq j_1 + j_2 + \dots +
	j_t$
	and $t=1,\dots,d$,
	\[
	f_{j_1+j_2+\dots+j_d+1}^{(s)}(x) = \bar{f}_{hat,u_{i_k}}(\bb_l^T x),
	\]
	and
	\[
	f^{(s)}_k(x)=1
	\]
	for $k=j_1+j_2+\dots+j_d+2, j_1+j_2+\dots+j_d+3, \dots,2^s$.
	As before,
	it is easy to see that
	$f_{net,k,j_1,\dots,j_d,\bb_l}$ is a neural network with $s+2$ hidden
	layers
	and at most
	\[
	6 \cdot 2^s, \,
	12 \cdot 2^s, \,
	2 \cdot 2^s,2^s,\,
	\dots,
	8,4
	\]
	neurons in the layers $1,2,\dots,s+2$, resp.
	Consequently, this
	network
	is contained in the class of all fully connected neural networks
	with $s+2$ hidden layers and $24 \cdot (N+1)$ neurons
	in each hidden layer.
	Furthermore it is easy to see that all weights are bounded in absolute
	value by $c_8 \cdot \max\{1, M/A, R^2\}$.

	\subsection{Definition of the output weights}
	\label{se3sub2}
	For given directions $\bb_l$ $(l=1,\dots,r)$
	we define our neural network estimate $\tilde{m}_n(x)$ by
	\[
	\tilde{m}_n(x)
	=
	\sum_{l=1,\dots,r}
	\sum_{k=1}^{M+1}
	\sum_{\substack{j_1, \dots, j_d \in \{0, \dots, N\}\\j_1 + \dots+j_d \leq N}}
	a_{i_k,j_1,\dots,j_d,\bb_l}
	\cdot
	f_{net,j_1,\dots,j_d,i_k,\bb_l}(x),
	\]
	where the coefficients $a_{i_k ,j_1,\dots,j_d,\bb_l}$ are chosen by minimizing
	\begin{equation}
          \label{se3sub2eq1}
	\frac{1}{n} \sum_{i=1}^n |Y_i-\tilde{m}_n(X_i)|^2
	+
	\frac{c_3}{n} \cdot
	\sum_{l=1}^r
	\sum_{k=1}^{M+1}
	\sum_{\substack{j_1, \dots, j_d \in \{0, \dots, N\}\\j_1 + \dots+j_d \leq N}}
	a_{i_k,j_1,\dots,j_d,\bb_l}^2
	\end{equation}
	for some constant $c_3>0$.
	This regularized linear least squares estimate can be computed by
	solving a linear equation system. To see this, set
	\[
	J= r \cdot (M+1) \cdot \left({N+d \atop d}  \right),
	\]
	let
	\begin{eqnarray*}
          &&
	  \{ B_{j} \, : \, j=1,\dots, J\}
          \\
          &&
	=
	\left\{
	f_{net,j_1,\dots,j_d,i_k,\bb_l}(x)
	\, : \,
	1 \leq l \leq r, \,
	1 \leq k \leq M+1
	\mbox{ and }
	0 \leq j_1 + \dots+j_d \leq N
	\right\}
	\end{eqnarray*}
	and set
	\[
	  {\mathbf B} 
          = (B_j(X_i))_{1 \leq i \leq n , 1 \leq j \leq J}
	\quad
	\mbox{and}
	\quad
	{\mathbf Y}=(Y_i)_{i=1,\dots,n}.
	\]
        As in Subsection \ref{se2sub3}
	it is easy to see
	that the vector of coefficients of our estimate 
	is the unique solution of the linear equation system
	\begin{equation}
	\label{se3eq9}
	\left(
	\frac{1}{n}
	{\mathbf B}^T  {\mathbf B} +
	\frac{c_3}{n} \cdot {\mathbf 1}
	\right)
	{\mathbf a}
	=
	\frac{1}{n}
	{\mathbf B}^T {\mathbf Y}.
	\end{equation}
	The value of (\ref{se3sub2eq1}) will be
	also
	less than or equal to
	the value which we get for coefficients equal to zero,
	hence we have
	\begin{equation}
	\label{se3eq9b}
	\frac{1}{n} ({\mathbf Y}- {\mathbf B a})^T ({\mathbf Y} - {\mathbf B
		a})
	+
	\frac{c_3}{n} \cdot{\mathbf a}^T {\mathbf a}
	\leq
	\frac{1}{n} \sum_{i=1}^n Y_i^2,
	\end{equation}
	which will allow us to derive a bound on the maximal absolute value
	of our coefficients.
	
	\subsection{Choice of the directions}
	\label{se3sub3}
	In order to choose $\bb_l$ $(l=1, \dots,r)$, we
	choose them $I_n$ times independent randomly according to
	a uniform distribution on $[-1,1]^d$,
	compute each time the corresponding outer weigths as in Subsection
	\ref{se3sub2}, and choose the directions and the
	corresponding outer weights for our estimate $\tilde{m}_n$,
	where the empirical $L_2$ risk of the estimate is minimal.
	
	\subsection{Rate of convergence}
	\label{se3sub4}
	
	\begin{theorem}
		\label{th2}
		Assume that the distribution of $(X,Y)$ satisfies
                (\ref{subgaus})
                for some constant $c_{4}>0$ and that the distribution of $X$ has bounded
		support $supp(X)$, and let $m(x)=\EXP\{Y|X=x\}$ be the corresponding
		regression function. Let $r \in \N$, $p>0$ and $C>0$,  and assume that the
		regression function satisfies
		\[
		m(x)=
		\sum_{l=1}^r g_l (\ba_l^T x)
		\quad (x \in \Rd)
		\]
		for some $(p,C)$--smooth functions
		$g_l:\R \rightarrow \R$
		and some $\ba_l \in \R^d$
		with $\|\ba_l\|=1$
		$(l=1, \dots, r)$.

		Define the estimate $\tilde{m}_n$ as in Subsections \ref{se3sub1}-\ref{se3sub3},
		where $\sigma$ is the logistic squasher and where
		$I_n=\lceil c_9\cdot (\log n)^2 \cdot n^{\frac{r \cdot d}{2p+1}} \rceil$, $N \geq p$,
		$M=M_n=\lceil c_{10} \cdot n^{1/(2p+1)} \rceil$,
		$R=R_n=n^{3}$
		and
		$A=A_n=(\log n)^{1/(6(N+d)}$. 
		Set
		$\beta_n = c_{6} \cdot \log(n)$ for some suitably large constant $c_{6}>0$
		and define $m_n$ by
		\[
		m_n(x)=T_{\beta_n} \tilde{m}_n(x).
		\]
		Then
		$m_n$ satisfies for
		$n$ suffciently large
		\[
		\mathbf E \int |m_n(x) - m(x)|^2 \PROB_X (dx) \leq
		c_{11} \cdot (\log n)^3 \cdot n^{-\frac{2p}{2p+1}},
		\]
		where $c_{11} >0$ does not depend on $n$.
	\end{theorem}

        \noindent
            {\bf Remark 2.} In order to compute our estimate, we
            have to solve $I_n$ times a linear equation system
            with a quadratic matrix of size $M_n$, for which
            computing time is proportional to
            \[
I_n \cdot M_n^2 \approx (\log n)^2 \cdot n^{\frac{r \cdot d+2}{2p+1}}.
            \]
            Hence in case
            \[
            r \cdot d < 4 \cdot p 
            \]
computing time is $O(n^2)$, so in case that
            the number $r$ of terms in the projection pursuit model
            and the dimension $d$ of $X$ are not too large,
            our estimate can be computed in  $O(n^2)$ time.

\section{Application to simulated data}
\label{se4}
            
In this section we illustrate the finite sample size performance of
our newly proposed
estimate by applying it
to simulated data
using the software \textit{MATLAB}.

The simulated data which we use is defined as follows: 
We choose $X$ uniformly distributed on $[-1,1]^d$, 
where $d$ is the dimension of the input,
$\epsilon$ standard normal and independent of $X$, 
and we define $Y$ by
\begin{equation}
\label{se4eq1}
Y=m_j(X) + \sigma \cdot \lambda_j \cdot \epsilon,
\end{equation}
where $m_j:[-1,1]^d \rightarrow \R$ is described below,
$\lambda_j>0$ is a scaling value defined below and $\sigma$
is chosen from $\{0.05,0.10\}$ $(j \in \{1,2,3,4\})$.
As regression functions we use 
\begin{eqnarray*}
%m_1(x_1,x_2,x_3,x_4)=2 \cdot \cos \left(
%\frac{\pi}{ x_1-0.5\cdot x_2 + 0.3 \cdot x_3+ 0.5 \cdot x_4  -3}
%\right),
m_1(x_1,x_2) 
&=&
 \log(0.2 \cdot x_1 + 0.9 \cdot x_2) + \cos \left( \frac{\pi}{\log (0.5\cdot x_1 + 0.3 \cdot x_2)} \right)
 \\
 &&
 \vspace{2mm}
 + \exp \left( \frac{1}{50} \cdot (0.7 \cdot x_1 + 0.7 \cdot x_2) \right)
 + \frac{\tan \left( \pi \cdot (0.1 \cdot x_1  + 0.3 \cdot x_2)^4 \right)}{(0.1\cdot x_1 + 0.3 \cdot x_2)^2},
\end{eqnarray*}
%so $m_1$ satisfies a single index model, and
%
\begin{eqnarray*}
	m_2(x_1,x_2,x_3,x_4)
	&=&
	\tan \left( \sin \left( \pi \cdot (0.2 \cdot x_1 + 0.5 \cdot x_2 - 0.6 \cdot x_3 +0.2 \cdot x_4) \right) \right)
	\\
	&&
	\vspace{2mm}
	+ \left( 0.5 \cdot \left( x_1 + x_2 + x_3 + x_4 \right) \right)^3
		\\
	&&
	\vspace{2mm}
	+ \frac{1}{\left( 0.5 \cdot x_1 + 0.3 \cdot x_2 - 0.3 \cdot x_3 + 0.25 \cdot x_4 \right)^2 + 4},
\end{eqnarray*}
%
% m_3 still running
%
\begin{eqnarray*}
	%actually m4
	m_3(x_1,x_2,x_3,x_4,x_5) 
	&=&
	\log \left( 0.5 \cdot (x_1 + 0.3 \cdot x_2 + 0.6 \cdot x_3 + x_4 - x_5)^2 \right)
	\\
	&&
	\vspace{2mm}
	+ \sin \left( \pi \cdot (0.7 \cdot x_1 + x_2 - 0.3 \cdot x_3 - 0.4 \cdot x_4 - 0.8 \cdot x_5) \right)
	\\
	&&
	\vspace{2mm}
	+ \cos \left( \frac{\pi}{1 + \sin(0.5 \cdot (x_2 + 0.9 \cdot x_3 - x_5))} \right)
\end{eqnarray*}
and
\begin{eqnarray*}
	%actually m5
&&
	m_4(x_1,x_2,x_3,x_4,x_5,x_6)
\\
	&&=
	\exp \left( 0.2 \cdot (x_1 + x_2 + x_3 + x_4 + x_5 + x_6) \right)
	\\
	&&
	\hspace*{0.4cm}
	+ \sin \left( \frac{\pi}{2} \cdot (x_1 - x_2 -x_3 + x_4 -x_5 - x_6) \right)
	\\
	&&
	\hspace*{0.4cm}
	+ \frac{1}{\left( 0.3 \cdot x_1 - 0.2 \cdot x_2 + 0.8 \cdot x_3 - 0.5 \cdot x_4 + 0.6 \cdot x_5 - 0.2 \cdot x_6 \right)^2 + 6}
	\\
	&&
	\hspace*{0.4cm}
	+ 0.5 \cdot (x_1 + x_3 - x_5)^3
\end{eqnarray*}
%
%hence $m_2$ satisfies a single index model with $r=2$ terms.
$\lambda_j$ is chosen approximately as IQR of a sample of size
$100$ of $m(X)$,
and we use the values 
%$\lambda_1=4.000$ and $\lambda_2=1.3584$.
$\lambda_1=5.04$,
$\lambda_2=5.57$,
$\lambda_3=6.8$,
and $\lambda_4=3.71$.
From this distribution we generate a sample of size $n=100$ and apply our
newly proposed neural network regression estimate and compare our
results to that of 
six alternative
regression estimates on the same data. Then we compute the $L_2$ errors
of these estimates approximately by using the empirical $L_2$
error
$\varepsilon_{L_2,\bar{N}}(\cdot)$
on an independent sample of $X$ of size $\bar{N}=10,000$.
Since this error strongly depends on the behavior of the correct function $m_j$, we consider it in relation to the error of the simplest estimate for $m_j$ we can think of, a completely constant function (whose value is the average of the observed data according to the least squares approach). Thus, the scaled error measure we use for evaluation of the estimates is $\varepsilon_{L_2,\bar{N}}(m_{n,i})/\bar{\varepsilon}_{L_2,\bar{N}}(avg)$, where $\bar{\varepsilon}_{L_2,\bar{N}}(avg)$ is the median of $50$ independent realizations of the value
one obtains if one plugs
the average of $n$ observations into $\varepsilon_{L_2,\bar{N}}(\cdot)$. To a certain extent, this quotient can be interpreted as the relative part of the error of the constant estimate that is still contained in the more sophisticated approaches.
The resulting scaled
errors of course depend on the random sample of $(X,Y)$, and to
be able to compare these values nevertheless we repeat the whole
computation $50$ times and report the median and the interquartile
range of the $50$ scaled errors for each of our estimates.

We choose the parameters for each of the estimates
by splitting of the sample. Here we split our sample in a learning
sample of size $n_l=0.8 \cdot n$ and a testing sample of size
$n_t=0.2 \cdot n$. We compute the estimate for all parameter
values from the sets described below using the learning sample, compute
the corresponding empirical $L_2$ risk on the testing sample
and choose the parameter value which leads to the minimal
empirical $L_2$ risk on the testing sample.

Our first three estimates are fully connected neural network estimates where the number of layers is fixed and the number of neurons per layer is chosen adaptively. 
The estimate {\it fc-neural-1} has one hidden layer, 
estimate {\it fc-neural-3} has three hidden layers, 
estimate {\it fc-neural-6} has six hidden layers and the number of
neurons 
per layer is chosen from the set 
$\{5, 10, 25, 50, 75\}$, 
$\{3,6,9,12,15\}$, 
$\{2,4,6,8,10\}$, respectively.

Our fourth estimate {\it kernel} 
is the Nadaraya-Watson kernel estimate 
with so-called naive kernel
where the bandwith is chosen from the set
$\{2^k : k \in \{-5,-4,\dots,5\}\}$.

Our fifth estimate {\it neighbor}
is a nearest neighbor estimate where the number
of nearest neighbors is chosen from the set $\{1,2,3\} \cup \{4,8,12,16, \dots, 4 \cdot \lfloor \frac{n_{l}}{4} \rfloor\}$.

Our sixth estimate {\it RBF} is the interpoland with radial basis functions
where the radial basis functions $\Phi(r) = (1 - r)_+^6 \cdot (35 \cdot r^2 + 18 \cdot r + 3)$ is used
and the scaling radius is chosen adaptively.

Our seventh estimate {\it MARS} is a method which makes use of multivariate adaptive regression splines. For this estimate we use the MATLAB ARESLab toolbox.

Our last estimate {\it proj-neural} is our newly proposed
neural network estimate presented in this paper.
Here the following parameters of the estimate are fixed: 
$N$ is set to 2, 
$A$ is set to 1, 
and $R$ is set to $10^6$, 
and $r$ is set to 4. %%% change
The parameter $M$ of the estimate is chosen
from the set $\{2,4,8,16\}$.  In order
to accelerate the computation of this estimate we use
only $I_n=50$ random choices for the vectors of directions
in the computation of the estimate $m_1$ with noise value $0.05$ 
and $I_n=400$ random choices for the vectors of directions
in the computation of each of the other estimates
for each parameter
value.

The results are summarized in Table \ref{se4ta1} and in Table
\ref{se4ta2}.
As we can see from the reported scaled errors, our
newly proposed neural network estimate outperforms
all other estimates in four out of eight cases. 
In the other settings our proposed neural network is only
outperformed by a fully connected network.
Still, in these scenarios we can see that the scaled error results of
our estimate are able to compete with those of the fully connected
neural networks, in the sense that the values of the former lie within
a small range of the best error value. 

\begin{table}[htbp]                                  
	\centering
	\makebox[\textwidth][c]{                                     
		\begin{tabular}{|c|c|c|c|c|}
			\hline
			&\multicolumn{2}{|c|}{$m_1$} & \multicolumn{2}{|c|}{$m_2$}\\
			\hline
			\textit{noise} & {$5\%$} & {$10\%$} & {$5\%$} & {$10\%$}\\
			\hline
			$\bar{\varepsilon}_{L_2,\bar{N}}(avg)$ & {$2.586$} & {$2.5892$} &  {$1.7504$} &  {$1.7504$} \\
			\hline                                         
			\textit{approach} & median (IQR) & median (IQR) & median (IQR) & median (IQR) \\                  
			\hline
			fc-neural-1 &  $0.0564$ $(0.035)$ & $0.0941$ $(0.036)$   & $0.0335$  $(0.033)$ & $\textbf{0.0816}$ $\textbf{(0.036)}$  \\
			fc-neural-3 & $0.0717$ $(0.0292)$ & $0.0898$ $(0.034)$  & $0.0373$  $(0.057)$ &  $0.1055$ $(0.054)$  \\
			fc-neural-6 & $0.0802$ $(0.041)$ & $0.0985$ $(0.044)$   & $0.0387$  $(0.018)$ &  $0.0914$ $(0.048)$  \\
			kernel & $0.1639$ $(0.052)$ & $0.1731$ $(0.056)$    & $0.1448$  $(0.058)$ &  $0.1631$ $(0.062)$  \\
			neighbor & $0.1254$ $(0.033)$ & $0.1452$ $(0.044)$   & $0.1449$  $(0.038)$ & $0.1824$ $(0.06)$ \\
			RBF & $0.3740$ $(0.413)$ & $0.9022$ $(1.085)$   & $0.1622$ $(0.118)$ & $0.43$ $(0.185)$    \\
			MARS & $0.0907$ $(0.08)$ & $0.1161$ $(0.105)$   & $0.0625$ $(0.084)$ & $0.1453$ $(0.106)$    \\
			proj-neural & $\textbf{0.047}$ $\textbf{(0.040)}$ & $\textbf{0.0732}$ $\textbf{(0.040)}$   & $\textbf{0.0318}$  $\textbf{(0.017)}$ & $0.1113$ $(0.097)$  \\
			\hline
		\end{tabular}
	}
	\caption{Median and IQR of the scaled empirical $L_2$ error of estimates for $m_1$ and $m_2$ for sample size $n=100$.}                
	\label{se4ta1}                            
\end{table}

\begin{table}[htbp]                                  
	\centering
	\makebox[\textwidth][c]{                                     
		\begin{tabular}{|c|c|c|c|c|}
			\hline
			&\multicolumn{2}{|c|}{$m_3$} & \multicolumn{2}{|c|}{$m_4$}\\
			\hline
			\textit{noise} & {$5\%$} & {$10\%$} & {$5\%$} & {$10\%$}\\
			\hline
			$\bar{\varepsilon}_{L_2,\bar{N}}(avg)$ & {$4.4779$} & {$4.4796$} &  {$0.5401$} &  {$0.54$} \\
			\hline                                         
			\textit{approach} & median (IQR) & median (IQR) & median (IQR) & median (IQR) \\                  
			\hline
			fc-neural-1 & $0.2938$ $(0.070)$ & $0.3431$ $(0.279)$   & $\textbf{0.1328}$  $\textbf{(0.067)}$ & $0.3172$ $(0.222)$  \\
			fc-neural-3 & $0.2451$ $(0.201)$ & $\textbf{0.2567}$ $\textbf{0.207}$  & $0.1787$  $(0.156)$ &  $0.3800$ $(0.399)$  \\
			fc-neural-6 & $\textbf{0.2141}$ $\textbf{(0.141)}$ & $0.2688$ $(0.207)$   & $0.1393$  $(0.125)$ &  $0.3272$ $(0.313)$  \\
			kernel & $0.3516$ $(0.070)$ & $0.3713$ $(0.077)$    & $0.8009$  $(0.095)$ &  $0.7961$ $(0.087)$  \\
			neighbor & $0.3406$ $(0.066)$ & $0.3602$ $(0.075)$   & $0.4502$  $(0.107)$ & $0.5109$ $(0.128)$ \\
			RBF & $0.3668$ $(0.234)$ & $0.5162$ $(0.112)$   & $0.2017$ $(0.080)$ & $0.6378$ $(0.312)$    \\
			MARS & $0.4300$ $(0.708)$ & $0.4259$ $(0.393)$   & $0.5875$ $(0.288)$ & $0.6816$ $(0.328)$    \\
			proj-neural & $0.236$ $(0.075)$ & $0.2822$ $(0.132)$   & $0.1433$  $(0.048)$ & $\textbf{0.2488}$ $\textbf{(0.171)}$  \\
			\hline
		\end{tabular}
	}
	\caption{Median and IQR of the scaled empirical $L_2$ error of estimates for $m_3$ and $m_4$ for sample size $n=100$.}                
	\label{se4ta2}                            
\end{table}

\section{Proofs}
\label{se5}

\subsection{Approximation results for neural networks}
We will use the following assumption on the activation function
of our neural network.
\begin{definition}
  \label{se5de1}
  Let $N \in \N_0$.
  A function $\sigma : \mathbb{R} \to [0, 1]$ is called \textbf{N-admissible},
  if it is nondecreasing and Lipschitz continuous and if, in addition,
  the following three conditions are satisfied:
\begin{itemize}
\item[(i)] The function $\sigma$ is $N+1$ times continuously differentiable with bounded derivatives.
\item[(ii)] A point $t_{\sigma} \in \mathbb{R}$ exists, where all
  derivatives up to  order $N$ of $\sigma$ are
 nonzero.
\item[(iii)] If $y > 0$, the relation $|\sigma(y) - 1| \leq \frac{1}{y} $ holds. If $y < 0$, the relation $|\sigma(y)| \leq \frac{1}{|y|}$ holds. \label{3}
\end{itemize}
\end{definition}
\noindent
It is easy to see that the logistic squasher (\ref{se2eq4})
is $N$--admissible for any $N \in \N$ (cf., e.g., Bauer and Kohler (2019)).

\begin{lemma}
  \label{le1}
Let $\sigma:\R \rightarrow \R$ be a function, let $R,a>0$.

\noindent
{\bf a)}
Assume that $\sigma$ is two times continuously differentiable and
let
$t_{\sigma,id} \in \R$ be such that $\sigma^\prime(t_{\sigma,id}) \neq
0$.
Then
\[
f_{id}(x)
=
\frac{R}{\sigma^\prime(t_{\sigma,id})}
\cdot
\left(
\sigma \left(
\frac{x}{R} + t_{\sigma,id}
\right)
-
\sigma(t_{\sigma,id})
\right)
\]
satisfies for any $x \in [-a,a]$:
\[
| f_{id}(x)-x|
\leq
\frac{
\| \sigma^{\prime \prime}\|_{\infty} \cdot a^2
}{
2 \cdot |\sigma^\prime(t_{\sigma,id})|
}
\cdot
\frac{1}{R}.
\]

\noindent
{\bf b)}
Assume that $\sigma$ is three times continuously differentiable and
let
$t_{\sigma,sq} \in \R$ be such that $\sigma^{\prime \prime}(t_{\sigma,sq}) \neq
0$.
Then
\[
f_{sq}(x)
=
\frac{R^2}{\sigma^{\prime \prime}(t_{\sigma,sq})}
\cdot
\left(
\sigma \left(
\frac{2x}{R} + t_{\sigma,sq}
\right)
-
2 \cdot
\sigma \left(
\frac{x}{R} + t_{\sigma,sq}
\right)
+
\sigma(t_{\sigma,sq})
\right)
\]
satisfies for any $x \in [-a,a]$:
\[
| f_{sq}(x)-x^2|
\leq
\frac{
5 \cdot
\| \sigma^{\prime \prime \prime}\|_{\infty} \cdot a^3
}{
3 \cdot |\sigma^{\prime \prime}(t_{\sigma,sq})|
}
\cdot
\frac{1}{R}.
\]
\end{lemma}

\noindent
\textbf{Proof.}
   The result follows in a straightforward way from the proof of
   Theorem 2 in Scarselli and Tsoi (1998), cf. Lemma 1 in Kohler,
   Krzy\.zak and Langer (2019). \hspace*{1cm} \hfill $\Box$

\noindent
    {\bf Remark 3.} In case of the logistic squasher it is easy to see
    that with
the choice $t_{\sigma,id}=0$ the network $f_{id}$ in
Lemma \ref{le1} is given by (\ref{se2eq6}).

\begin{lemma}
\label{le2}
Let $\sigma: \mathbb{R} \to [0,1]$ be 2-admissible according to
Definition \ref{se5de1}. Then for any $R>0$ and any $a >0$ the neural network
\begin{eqnarray*}
  f_{mult}(x,y) &=&
\frac{R^2}{4 \cdot \sigma^{\prime \prime}(t_\sigma)}
\cdot
\Bigg(
\sigma \left(
\frac{2 \cdot (x+y)}{R} + t_\sigma
\right)
-
2 \cdot
\sigma \left(
\frac{x+y}{R} + t_\sigma
\right)
\\
&&
\hspace*{2cm}
-
\sigma \left(
\frac{2 \cdot (x-y)}{R} + t_\sigma
\right)
+
2 \cdot
\sigma \left(
\frac{x-y}{R} + t_\sigma
\right)
\Bigg)
\end{eqnarray*}
satisfies for any $x \in [-a,a]$:
\begin{equation*}
|f_{mult}(x,y) - x \cdot y| \leq
\frac{20 \cdot \|
\sigma^{\prime \prime \prime}
\|_{\infty} \cdot a^3}{
3 \cdot |\sigma^{ \prime \prime} (t_\sigma)|}
 \cdot \frac{1}{R}.
\end{equation*}
\end{lemma}

\noindent {\bf Proof.}
See Lemma 2 in Kohler, Krzy\.zak and Langer (2019).
 \hfill $\Box$

\noindent
    {\bf Remark 4.} In case of the logistic squasher it is easy to see
    that with
the choice $t_{\sigma}=1$ the network $f_{mult}$ in
Lemma \ref{le2} is given by (\ref{se2eq7}).

\begin{lemma}\label{le3}
  Let $\sigma: \mathbb{R} \to [0,1]$ be 2-admissible
  according to Definition \ref{se5de1}.
Let $f_{mult}$ be the neural network from Lemma \ref{le2}
    and let $f_{id}$ be the network from Lemma \ref{le1}.
Assume
\begin{equation}
\label{le3eq1}
a \geq 1 \quad \mbox{and} \quad
R
\geq
\frac{
\| \sigma^{\prime \prime}\|_{\infty}\cdot a
}{
2 \cdot | \sigma^\prime (t_{\sigma.id})|
}.
\end{equation}
Then the neural network
\begin{eqnarray*}
  f_{ReLU}(x) &=&
f_{mult} \left(f_{id}(x),\sigma \left(R \cdot x \right) \right)
  \\
  &=&
  \sum_{k=1}^4 d_k \cdot \sigma \left(\sum_{i=1}^2 b_{k,i} \cdot \sigma(a_i \cdot x + t_{\sigma})
  +
  b_{k,3} \cdot \sigma (a_3 \cdot x)
  +  t_{\sigma} \right)
\end{eqnarray*}
satisfies
\begin{equation*}
|f_{ReLU}(x) - \max\{x,0\}| \leq56 \cdot \frac{
\max \left\{
\| \sigma^{\prime \prime}\|_{\infty},
\| \sigma^{\prime \prime \prime}\|_{\infty},1
\right\}
}{
\min \left\{
2 \cdot |\sigma^\prime (t_{\sigma.id})|,
|\sigma^{\prime \prime} (t_{\sigma})|, 1
\right\}
} \cdot a^3 \cdot \frac{1}{R}
\end{equation*}
for all $x \in [-a, a]$.
\end{lemma}

\noindent {\bf Proof.}
See Lemma 3 in Kohler, Krzy\.zak and Langer (2019).\hfill $\Box$
 
\begin{lemma}\label{le4}
  Let $M \in \N$ and
  let $\sigma: \mathbb{R} \to [0,1]$ be 2-admissible  according to
  Definition \ref{se5de1}. Let $a>0$ and
\[
R
\geq
\frac{
\| \sigma^{\prime \prime}\|_{\infty}\cdot (M+1)
}{
2 \cdot | \sigma^\prime (t_{\sigma.id})|
},
\]
let $y \in [-a,a]$
and let $f_{ReLU}$ be the neural network of Lemma \ref{le3}.
Then the network
\begin{eqnarray*}
f_{hat,y}(x)
&=&
f_{ReLU} \left(
\frac{M}{2a} \cdot (x-y) + 1
\right)
-
2
\cdot
f_{ReLU} \left(
\frac{M}{2a} \cdot (x-y)
\right)
\\
&&
+
f_{ReLU} \left(
\frac{M}{2a} \cdot (x-y) -1
\right)
\end{eqnarray*}
satisfies
\[
\left| f_{hat,y}(x)-
\left(
1 - \frac{M}{2a} \cdot |x-y|
\right)_+
\right|
\leq
1792 \cdot \frac{
\max \left\{
\| \sigma^{\prime \prime}\|_{\infty},
\| \sigma^{\prime \prime \prime}\|_{\infty},1
\right\}
}{
\min \left\{
2 \cdot |\sigma^\prime (t_{\sigma.id})|,
|\sigma^{\prime \prime} (t_{\sigma})|, 1
\right\}
} \cdot M^3 \cdot \frac{1}{R}
\]
for all $x \in [-a, a]$.
\end{lemma}

\noindent {\bf Proof.}
  Since
  \[
  (1-\frac{M}{2a} \cdot|x|)_+=\max\{\frac{M}{2a} \cdot x+1,0\}-2 \cdot \max\{\frac{M}{2a} \cdot
  x,0\}+\max\{\frac{M}{2a} \cdot x-1,0\}
  \quad (x \in \R)
  \]
  the result is an easy consequence of Lemma \ref{le3} (applied with
  $M+1$
instead of $a$).
\hfill $\Box$

\begin{lemma}
\label{le5}
Let $M \in \N$ and
  let $\sigma: \mathbb{R} \to [0,1]$ be 2-admissible  according to
  Definition \ref{se5de1}. Let $a \geq 1$ and
\begin{eqnarray}
\label{le5eq1}
R
& \geq&
\max\Bigg\{
\frac{
\| \sigma^{\prime \prime}\|_{\infty}
\cdot  (M+1)
}{
2 \cdot  | \sigma^\prime (t_{\sigma,id})|
}
,
\frac{
9 \cdot
\| \sigma^{\prime \prime}\|_{\infty}\cdot a
}{
 | \sigma^\prime (t_{\sigma,id})|
}
,\\
&&
\frac{20 \cdot \| \sigma^{\prime \prime \prime}\|_\infty}{
3 \cdot |\sigma^{\prime\prime}(t_{\sigma})|}
\cdot 3^{3 \cdot 3^s} \cdot a^{3 \cdot 2^s},
1792 \cdot \frac{
\max \left\{
\| \sigma^{\prime \prime}\|_{\infty},
\| \sigma^{\prime \prime \prime}\|_{\infty},1
\right\}
}{
\min \left\{
2 \cdot |\sigma^\prime (t_{\sigma,id})|,
|\sigma^{\prime \prime} (t_{\sigma})|, 1
\right\}
} \cdot M^3 
\Bigg\}
\nonumber
\end{eqnarray}
and let $y \in [-a,a]^d$. Let $N \in \N$ and let
$j_1, \dots, j_d \in \N_0$ such that $j_1+\dots+j_d \leq N$, and set
$s=\lceil \log_2(N+d) \rceil$.
Let $f_{id}$, $f_{mult}$ and $f_{hat,z}$ (for $z \in \R$) be the
  neural
networks defined in Lemma \ref{le1}, Lemma \ref{le2} and Lemma
\ref{le4}, resp.
Define the network $f_{net,j_1,\dots,j_d,y}$ by
\[
f_{net,j_1,\dots,j_d,y}(x)=f_1^{(0)}(x),
\]
where $f_1^{(0)}$ is defined by backward recursion as follows:
\[
f_k^{(l)}(x)
=
f_{mult}
\left(
f_{2k-1}^{(l+1)}(x),f_{2k}^{(l+1)}(x)
\right)
\]
for $k \in \{1, 2, \dots, 2^l\}$ and $l \in \{0,\dots,s-1\}$,
and
\[
f_k^{(s)}(x)=f_{id}(f_{id}(x^{(l)}-y^{(l)}))
\]
for $j_1 + j_2 + \dots + j_{l-1}+1 \leq k \leq j_1 + j_2 + \dots +
j_l$
and $l=1,\dots,d$,
\[
f_{j_1+j_2+\dots+j_d+k}^{(s)}(x) = f_{hat,y^{(k)}}(x^{(k)})
\]
for $k=1, \dots, d$, and
\[
f^{(s)}_k(x)=1
\]
for $k=j_1+j_2+\dots+j_d+d+1, j_1+j_2+\dots+j_d+d+2, \dots,2^s$.
Then we have for any $x \in [-a,a]^d$:
\begin{eqnarray*}
&&
\left|
f_{net,y}(x)
-
(x^{(1)}-y^{(1)})^{j_1} \cdots (x^{(d)}-y^{(d)})^{j_d}
\prod_{j=1}^d(1 - \frac{M}{2a} \cdot |x^{(j)} - y^{(j)}|)_+
\right|
\\
&&
\leq
c_{12} \cdot 3^{3 \cdot 3^s} \cdot a^{3 \cdot 2^s} \cdot M^3 \cdot \frac{1}{R}.
\end{eqnarray*}
\end{lemma}

\noindent
    {\bf Proof.}
    The result follows from Lemma \ref{le1}, Lemma \ref{le2}
    and Lemma \ref{le4} in a straightforward but technical way
    using an induction. A complete proof
    can be found in the Supplement.
\hspace*{2cm} \hfill $\Box$

\noindent
{\bf Remark 5.}
The result can be analogously stated for our estimate in the context of the projection pursuit model. The corresponding statement and a complete proof can be  found in the Supplement.

\subsection{Approximation of a projection pursuit model by piecewise
	polynomials}

\begin{lemma}
	\label{le1b}
	Let $p=q+s$ for some $q \in \N_0$ and $s \in (0,1]$. 
        Let $C>0$, $r \in \N$, $g_l:\R \rightarrow \R$
	$(p,C)$-smooth functions $(l=1, \dots, r)$ and $\ba_l \in \R^d$
	$(l=1, \dots, r)$.
	Set
	\[
	m(x)=
	\sum_{l=1}^r g_l (\ba_l^T x)
	\quad (x \in \Rd).
	\]
	For $\bb_l \in \R^d$ $(l=1, \dots,r)$ set
	\[
	g(x)=
	\sum_{l=1}^r \sum_{j=0}^q \frac{
		g_l^{(j)}(\bb_l^Tx)
	}{j!}
	\cdot
	((\ba_l-\bb_l)^T x)^j,
	\]
        where $g_l^{(j)}$ denotes the $j$--th derivative of $g_l$.
	Then we have for any $x \in \Rd$
	\[
	|m(x)-g(x)|
	\leq
	\frac{r \cdot d^{p} \cdot C}{q!} \cdot \|x\|_\infty^p \cdot \|\ba_l-\bb_l\|_\infty^p.
	\]
\end{lemma}

\noindent
{\bf Proof.}
By the proof of Lemma 11.1 in Gy\"orfi et al. (2002) we have for
any $z \in \R$
\[
\left|
g_l(u)
-
\sum_{j=0}^q
\frac{g_l^{(j)}(z)}{j!} \cdot (u-z)^j
\right|
\leq
\frac{1}{q!} \cdot C \cdot |u-z|^p
\quad (u \in \R).
\]
Applying this with $u=\ba^T_l x$ and $z=\bb_l^T x$ we get
\begin{eqnarray*}
	&&
	|m(x)-g(x)|
	\\
	&&
	\leq
	\sum_{l=1}^r
	\left|
	g_l (\ba_l^T x)
	-
	\sum_{j=0}^q
	\frac{g_l^{(j)}(\bb_l^T x)}{j!} \cdot (\ba_l^T x -\bb_l^T x)^j
	\right|
	\\
	&&
	\leq
	\sum_{l=1}^r
	\frac{1}{q!} \cdot C \cdot |\ba_l^T x - \bb_l^T x|^p
	\\
	&&
	\leq
	\frac{r \cdot d^{p} \cdot C}{q!} \cdot \|x\|_\infty^p \cdot \|\ba_l-\bb_l\|_\infty^p.
\end{eqnarray*}
\quad \hfill $\Box$

\begin{lemma}
	\label{le2b}
	Let $p=q+s$ for some $q \in \N_0$ and $s \in (0,1]$.
	Let $C>0$, $r \in \N$, $g_l:\R \rightarrow \R$
	$(p,C)$-smooth functions $(l=1, \dots, r)$ and $\ba_l\in \R^d$
	with $\|\ba_l\| = 1$ and $\bb_l  \in [-1,1]^d$
	$(l=1, \dots, r)$.
	Let $A \geq 1$, $M \in \N$, set
	\[
	u_i=- \sqrt{d} \cdot A + i \cdot \frac{2 \cdot \sqrt{d} \cdot A}{M}
	\quad (i=0, \dots, M)
	\]
	and $\{i_1, \dots, i_{M+1}\} = \{0, \dots, M\}$. Then there exist polynomials $p_{i_k,l}:\Rd \rightarrow \R$
	of total degree $q$, which depend on $\ba_l$ and $\bb_l$
	and where all coefficients are bounded in absolute
	value by
	\[
	 (q+1) \cdot 2^p \cdot d^{3p/2} \cdot A^p \cdot \max_{
		l \in \{1, \dots,r\}, j \in \{0, \dots, q\}}
	\|g_l^{(j)}\|_\infty,
	\cdot
	(\max\{ \max_{l=1, \dots,r }\|\ba_l-\bb_l\|_\infty,1\})^p, 
	\]
	such that we have for all $x \in [-A,A]^d$ 
	\begin{eqnarray*}
		&&
		\left|
		\sum_{l=1}^r \sum_{j=0}^q \frac{
			g_l^{(j)}(\bb_l^Tx)
		}{j!}
		\cdot
		((\ba_l-\bb_l)^T x)^j
		-
		\sum_{l=1}^r
		\sum_{k=1}^{M+1} p_{i_k,l}(x) \cdot \left(
		1 -  \frac{M}{2 \cdot \sqrt{d} \cdot A} \cdot |\bb_l^Tx - u_{i_k}|
		\right)_+
		\right|
		\\
		&&
		\leq
		r \cdot 2^{p} \cdot (p+1) \cdot C \cdot d^{3p/2} \cdot A^{2p} \cdot
		\left(
		\max \left\{
		\frac{1}{M}, \max_{l=1, \dots,r }\|\ba_l-\bb_l\|_\infty
		\right\}
		\right)^p.
	\end{eqnarray*}
\end{lemma}

\noindent
{\bf Proof.}
Let $p_{l,j,i_k}$ be the Taylor polynomial
of $g_l^{(j)}$ of degree $q-j$ around $u_{i_k}$.
Because of the $(p-j,C)$-smoothness of
$g_l^{(j)}$ Lemma 11.1 in Gy\"orfi et al. (2002)
implies that we have
\[
\left|
g_l^{(j)}(\bb_l^Tx)
-
p_{l,j,i_k}(\bb_l^Tx)
\right|
\leq
\frac{1}{(q-j)!} \cdot C \cdot |\bb_l^Tx-u_{i_k}|^{(p-j)}.  
\]
From this we can conclude for $x \in [-A,A]^d$
\begin{eqnarray*}
	&&
	\left|
	\frac{g_l^{(j)}(\bb_l^Tx)
	}{j!}
	\cdot
	((\ba_l-\bb_l)^T x)^j
	-
	\frac{p_{l,j,i_k}(\bb_l^Tx)
	}{j!}
	\cdot
	((\ba_l-\bb_l)^T x)^j
	\right|
	\\
	&&
	\leq
	\frac{1}{(q-j)!} \cdot C \cdot d^{j} \cdot A^j
	\cdot
	\left(
	\max\left\{
	|\bb_l^Tx-u_{i_k}|
	,
	\|\ba_l-\bb_l\|_\infty
	\right\}
	\right)^p.
\end{eqnarray*}
Using
\[
\sum_{k=1}^{M+1}
\left(
1 -  \frac{M}{2 \cdot \sqrt{d} \cdot A} \cdot |\bb_l^Tx - u_{i_k}|
\right)_+    
=1            
\]
for $x \in [-A,A]^d$,
this in turn implies for $x \in [-A,A]^d$
\begin{eqnarray*}
	&&
	\Bigg|
	\sum_{j=0}^q
	\frac{g_l^{(j)}(\bb_l^Tx)
	}{j!}
	\cdot
	((\ba_l-\bb_l)^T x)^j
	\\
	&&
	\hspace*{2cm}
	-
	\sum_{j=0}^q \sum_{k=1}^{M+1}
	\frac{p_{l,j,i_k}(\bb_l^Tx)
	}{j!}
	\cdot
	((\ba_l-\bb_l)^T x)^j
	\cdot
	\left(
	1 -  \frac{M}{2 \cdot \sqrt{d} \cdot A} \cdot |\bb_l^Tx - u_{i_k}|
	\right)_+    
	\Bigg|
	\\
	&&
	\leq
	\sum_{j=0}^q
	\sum_{k=1}^{M+1}
	\left|
	\frac{g_l^{(j)}(\bb_l^Tx)
	}{j!}     \cdot
	((\ba_l-\bb_l)^T x)^j
	-
	\frac{p_{l,j,i_k}(\bb_l^Tx)
	}{j!}
	\cdot
	((\ba_l-\bb_l)^T x)^j
	\right|
	\\
	&&
	\hspace*{3cm}
	\cdot
	\left(
	1 -  \frac{M}{2 \cdot \sqrt{d} \cdot A} \cdot |\bb_l^Tx - u_{i_k}|
	\right)_+
	\\
	&&
	\leq
	\sum_{j=0}^q
	\max_{ i_k \in \{0, \dots, M\}, \atop |b_l^Tx-u_{i_k}|
		\leq 2 \cdot \sqrt{d} \cdot A / M}
	\left|
	\frac{g_l^{(j)}(\bb_l^Tx)
	}{j!}     \cdot
	((\ba_l-\bb_l)^T x)^j
	-
	\frac{p_{l,j,i_k}(\bb_l^Tx)
	}{j!}
	\cdot
	((\ba_l-\bb_l)^T x)^j
	\right|
	\\
	&&
	\leq
        (q+1) \cdot C \cdot d^{q} \cdot A^q
	\left(
	\max \left\{
	\frac{2 \cdot \sqrt{d} \cdot A}{M}, \|\ba_l -\bb_l \|_\infty
	\right\}
	\right)^p.
\end{eqnarray*}
With
\[
p_{i_k,l}(x)
=
\sum_{j=0}^q
\frac{p_{l,j,i_k}(\bb_l^Tx)
}{j!}
\cdot
((\ba_l-\bb_l)^T x)^j            
\]
we get the assertion.
\quad \hfill $\Box$

\subsection{Auxiliary results}
\label{se5sub2}
\begin{lemma}
\label{le6}
Let $\beta_n = c_{6} \cdot \log(n)$ for some suitably large
constant $c_{6}>0$. Assume that the distribution of $(X,Y)$ satisfies
(\ref{subgaus})
for some constant $c_{4}>0$ and that the regression function $m$ is
bounded in absolute value. Let
$\F_n$ be a set of functions $f:\Rd \rightarrow \R$ and assume that
the estimate $m_n$ satisfies
\[
m_n=T_{\beta_n}\tilde{m}_n
\]
and
\[
\tilde{m}_n(\cdot)
=
\tilde{m}_n(\cdot,(X_1,Y_1),\dots,(X_n,Y_n))
\in \F_n
\]
and
\[
	\frac{1}{n}
	\sum_{i=1}^n |Y_i-\tilde{m}_n(X_i)|^2
	\leq
	\min_{l \in \Theta_n}
	\left(
	\frac{1}{n}
	\sum_{i=1}^n |Y_i-g_{n,l}(X_i)|^2+pen_n(g_{n,l}) 
	\right)
\]
for some nonempty parameter set $\Theta_n$,
some random functions $g_{n,l}: \R^d \rightarrow \R$
and some
deterministic penalty terms $pen_n(g_{n,l}) \geq 0$, where
the random function $g_{n,l}: \R^d \rightarrow \R$ depend
only on random variables
\[
\bb_1^{(1)}, \dots, \bb_r^{(1)}, \dots,
\bb_1^{(I_n)}, \dots, \bb_r^{(I_n)},
\]
which are independent of $(X_1,Y_1)$, $(X_2,Y_2)$, \dots
 
Then $m_n$ satisfies
\begin{align*}
  &  \mathbf E \int |m_n(x) - m(x)|^2 \PROB_X (dx) \\
  &\leq
                                                  \frac{c_{13}\cdot
                                                  (\log n)^2\cdot
    \left(
    \log\left(
    \sup_{x_1^n \in (supp(X))^n}
\mathcal{N}_1 \left(\frac{1}{n\cdot\beta_n}, \mathcal{F}_n, x_1^n\right)
\right)
+1
\right)
  }{n}\\
      &\quad + 2 \cdot  \EXP \left( \min_{l \in \Theta_n}
                                                  \int |g_{n,l}(x)-m(x)|^2 \PROB_X (dx) + pen_n(g_{n,l}) \right) 
\end{align*}
for $n>1$ and some constant $c_{13}>0$, which does not depend on $n$.
\end{lemma}
\noindent
    {\bf Proof.}
    This lemma follows in a straightforward way from the proof of
    Theorem 1 in Bagirov et al. (2009). A complete version of the proof
    is given in the Supplement.
    \hfill $\Box$\\

\noindent\newline
In order to bound the covering number
$\mathcal{N}_1 \left(\frac{1}{n\cdot\beta_n}, \mathcal{F}_n,
  x_1^n \right)$
we will use the following lemma.

\begin{lemma}
\label{le8}
Let $a >0$ and
let $d, N, J_n \in \N$ be such that $J_n \leq n^{c_{14}}$
and set  $\beta_n=c_6 \cdot \log n$.
Let $\sigma$ be $2$--admissible according to Definition \ref{se5de1}.
Let $\F$ be the set of all functions
defined by (\ref{se2eq1}), (\ref{se2eq2}) and (\ref{se2eq3}) where
$k_1=k_2=\dots=k_{L}=24 \cdot (N+d)$
and the weights are bounded in absolute value by $c_{15} \cdot n^{c_{16}}$.  Set
\[
\F^{(J_n)}=
\left\{
\sum_{j=1}^{J_n}
a_j \cdot f_j \, : \,
f_j \in \F \quad \mbox{and} \quad
\sum_{j=1}^{J_n} a_j^2 \leq c_{17} \cdot n^{c_{18}} 
\right\}.
\]
Then we have for $n>1$
\[
\log\left(
supp_{x_1^n \in [-a,a]^{d \cdot n}}
\mathcal{N}_1 \left(\frac{1}{n\cdot\beta_n}, \mathcal{F}^{(J_n)},
 x_1^n \right) \right)
\leq
c_{19} \cdot \log n \cdot J_n
\]
for some constant $c_{19}$ which depends only on $L$, $N$, $a$ and $d$.
\end{lemma}

\noindent
    {\bf Proof.} 
    Since the networks in $\F^{(J_n)}$ are linear combinations
    of $J_n$
    fully connected neural networks
    with $L$ hidden layers, a bounded number of neurons in each
    hidden layers and all weights bounded by a polynomial in $n$,
    the result follows by combining Lemma 16.6 in Gy\"orfi et al. (2002)
with
    Lemma 7 in the Supplement of Bauer et al. (2019).
    \quad \hfill $\Box$

\subsection{Proof of Theorem \ref{th1}}
Since $supp(\PROB_X)$ is bounded and $m$ is $(p,C)$--smooth,
we conclude that $m$ is bounded in absolute value, and we can assume
without loss of generality that $supp(X) \subseteq [-a_n,a_n]^d$
and $\|m\|_\infty \leq \beta_n$.

Let $\F$ be the set of all functions
defined by (\ref{se2eq1}), (\ref{se2eq2}) and (\ref{se2eq3}) where
$L=s+2=\lceil \log_2(N+d) \rceil+2$, where
$k_1=k_2=\dots=k_{L}=24 \cdot (N+d)$
and where the weights are bounded in absolute value by $n^{c_{20}}$.  Set
\[
\F^{(J_n)}=
\left\{
\sum_{j=1}^{J_n}
a_j \cdot f_j \, : \,
f_j \in \F \quad \mbox{and} \quad
\sum_{j=1}^{J_n} a_j^2 \leq c_{21} \cdot n
\right\}
\]
for $c_{21}$ chosen below,
where
\[
J_n = (M_n+1)^d
\cdot
|\{
(j_1,\dots,j_d) \, : \, j_1, \dots, j_d \in \{0, \dots, N\}, \, j_1 + \dots+j_d \leq N
\}|
\]
Then $J_n \leq (M_n+1)^d \cdot (N+1)^d$.

Let
\[
g_n(x)=\sum_{k=1}^{(M_n+1)^d}
 \sum_{\substack{j_1, \dots, j_d \in \{0, \dots, q\}\\j_1 + \dots+j_d
     \leq q}} \frac{1}{j_1! \cdots j_d!} \cdot
 \frac{\partial^{j_1+\dots+j_d}m}{\partial^{j_1}x^{(1)}
   \cdots \partial^{j_d}x^{(d)}}(x_{\bi_k})
\cdot f_{net,j_1,\dots,j_d,\bi_k}(x).
\]
Because of the $(p,C)$--smoothness of $m$ we know that
\begin{equation}
\label{pth1eq*}
\max_{
k \in \{1, \dots, (M_n+1)^d,
j_1, \dots, j_d \in \{0, \dots, q\}, j_1 + \dots+j_d
     \leq q
}
\left|
\frac{\partial^{j_1+\dots+j_d}m}{\partial^{j_1}x^{(1)}
   \cdots \partial^{j_d}x^{(d)}}(x_{\bi_k})
\right|
< \infty.
\end{equation}
Set
\begin{eqnarray}
\label{pth1eq2} \nonumber
c_{21} &=&
\max \Bigg\{  \frac{1 + \EXP\{Y^2\}}{c_3},
(N+1)^d \cdot
\max \Bigg\{
\left|
\frac{1}{j_1! \cdots j_d!} \cdot
 \frac{\partial^{j_1+\dots+j_d}m}{\partial^{j_1}x^{(1)}
   \cdots \partial^{j_d}x^{(d)}}(x_{\bi_k})
\right|^2
\quad :
\\
&&
\hspace*{4.3cm}
j_1, \dots, j_d \in \{0, \dots, q\}, \, j_1 + \dots+j_d
     \leq q
\Bigg\}
\Bigg\}
\end{eqnarray}
and let $A_n$ be the event that
\begin{equation}
\label{pth1eq3}
\frac{1}{n} \sum_{i=1}^n Y_i^2
\leq
1 + \EXP\{Y^2\}
\end{equation}
holds. Then
\[
\PROB(A_n^c)
\leq
\frac{\Var\{Y^2\}}{n}
\leq
\frac{c_{22}}{n}
\]
by Chebychev inequality.

Set $\hat{m}_n=T_{\beta_n} \tilde{m}_n =m_n$ in case that $A_n$ holds and set
$\hat{m}_n=T_{\beta_n} g_n$ otherwise. Then
\begin{eqnarray*}
&&
 \mathbf E \int |m_n(x) - m(x)|^2 \PROB_X (dx)
\\
&&
\leq
4 \beta_n^2 \cdot \PROB\{A_n^c\}
+
 \mathbf E \{ \int |m_n(x) - m(x)|^2 \PROB_X (dx) \cdot 1_{A_n} \}
\\
&&
\leq
\frac{4 \cdot c_{22} \cdot \beta_n^2}{n}
+
 \mathbf E \int |\hat{m}_n(x) - m(x)|^2 \PROB_X (dx)
.
\end{eqnarray*}
The definition of the estimate $\tilde{m}_n$ implies
\[
\tilde{m}_n(x) = \sum_{j=1}^{J_n} \hat{a}_j \cdot f_j
\]
for some $f_j \in \F$ and some $\hat{a}_j$ satisfying
\[
\sum_{j=1}^{J_n} \hat{a}_j^2 \leq \frac{1}{n} \sum_{i=1}^n Y_i^2 \cdot
\frac{n}{c_3}.
\]
Hence on $A_n$ we have
\[
\sum_{j=1}^{J_n} \hat{a}_j^2 \leq \frac{1 + \EXP Y^2}{c_3} \cdot n,
\]
and consequently  we can assume w.l.o.g. that
$m_n$ satisfies $m_n=T_{\beta_n} \bar{m}_n$ for some
$\bar{m}_n \in \F^{(J_n)}$. And since
\begin{eqnarray*}
&&
 \frac{1}{n}
\sum_{i=1}^n |Y_i-\tilde{m}_n(X_i)|^2
\\
&&
\leq
 \frac{1}{n}
\sum_{i=1}^n |Y_i-\tilde{m}_n(X_i)|^2
+
\frac{c_3}{n} \cdot
\sum_{k=1}^{(M_n+1)^d}
\sum_{\substack{j_1, \dots, j_d \in \{0, \dots, N\}\\j_1 + \dots+j_d \leq N}}
a_{\bi_k,j_1,\dots,j_d}^2
\\
&&
\leq
\frac{1}{n}
\sum_{i=1}^n |Y_i-g_n(X_i)|^2
\\
&&
\hspace*{0.5cm}
+
\frac{c_3}{n} \cdot
\sum_{k=1}^{(M_n+1)^d}
\sum_{\substack{j_1, \dots, j_d \in \{0, \dots, q\}\\j_1 + \dots+j_d \leq q}}
\left|
\frac{1}{j_1! \cdots j_d!} \cdot
 \frac{\partial^{j_1+\dots+j_d}m}{\partial^{j_1}x^{(1)}
   \cdots \partial^{j_d}x^{(d)}}(x_{\bi_k})
\right|^2
\end{eqnarray*}
(by definition of $\tilde{m}_n$) and (\ref{pth1eq*}), we also have
\[
 \frac{1}{n}
\sum_{i=1}^n |Y_i-\bar{m}_n(X_i)|^2
\leq
 \frac{1}{n}
\sum_{i=1}^n |Y_i-g_n(X_i)|^2
+ c_{23} \cdot \frac{(M_n+1)^d}{n}.
\]
Set
\begin{eqnarray*}
P_n(x)&=&\sum_{k=1}^{(M_n+1)^d}
 \sum_{\substack{j_1, \dots, j_d \in \{0, \dots, q\}\\j_1 + \dots+j_d
     \leq q}} \frac{1}{j_1! \cdots j_d!} \cdot
 \frac{\partial^{j_1+\dots+j_d}m}{\partial^{j_1}x^{(1)}
   \cdots \partial^{j_d}x^{(d)}}(x_{\bi_k})
\\
&&
\hspace*{1cm}
\cdot
(x^{(1)}-x_{\bi_k}^{(1)})^{j_1} \cdots (x^{(d)}-x_{\bi_k}^{(d)})^{j_d}
\cdot
 \prod_{j=1}^d(1- \frac{M_n}{2a} \cdot |x^{(j)} - x_{\bi_k}^{(j)}|)_+.
\end{eqnarray*}

Application of Lemma \ref{le6} (with $|\Theta_n|=1$ and $g_{n,1}=g_n$
deterministic)
yields
\begin{align*}
  &  \mathbf E \int |m_n(x) - m(x)|^2 \PROB_X (dx)\\
  &\leq \frac{c_{23}\cdot
                                                 (\log n)^2\cdot
\left(
\log\left(
\sup_{x_1^n \in supp(X)^n}
\mathcal{N}_1 \left(\frac{1}{n\cdot\beta_n}, \mathcal{F}^{(J_n)}, x_1^n \right)
\right)
+1
\right)
  }{n}\\
  &\quad + 2 \cdot \int |g_n(x)-m(x)|^2 \PROB_X (dx)
  +
  2 \cdot c_{21} \cdot \frac{(M_n+1)^d}{n}.
\end{align*}
By Lemma \ref{le8} we know that
\begin{eqnarray*}
  &&
  \frac{ c_{23} \cdot
  \log(n)^2\cdot  \left(\log\left(
  \sup_{x_1^n \in supp(X)^n}
\mathcal{N}_1 \left(\frac{1}{n\cdot\beta_n}, \mathcal{F}^{(J_n)}, x_1^n \right)
\right) +1 \right)}{n}
\\
&&
\leq
c_{24} \cdot
\frac{(\log n)^3\cdot (N+1)^d \cdot (M_n+1)^d}{n}.
\end{eqnarray*}
Furthermore we have
\[
\int |g_n(x)-m(x)|^2 \PROB_X (dx)
\leq
2 \cdot \sup_{x \in [-a_n,a_n]^d}
|g_n(x)-P_n(x)|^2
+
2 \cdot \sup_{x \in [-a_n,a_n]^d}
|P_n(x)-m(x)|^2.
\]
By Lemma \ref{le5} we know
\begin{eqnarray*}
 \sup_{x \in [-a_n,a_n]^d}
|g_n(x)-P_n(x)|
&\leq&
(M_n+1)^d
\cdot
(q+1)^d
\cdot c_{25} \cdot 
a_n^{6(N+d)} \cdot M_n^3
\frac{1}{R_n}
\\
&
\leq & (M_n+1)^d \cdot (q+1)^d
\cdot c_{25} 
\cdot (\log n) \cdot \frac{ M_n^3}{R_n},
%\leq const \cdot \frac{\log n}{n}
\end{eqnarray*}
and Lemma 5 in
 Schmidt--Hieber (2019)
implies
\[
 \sup_{x \in [-a_n,a_n]^d}
|P_n(x)-m(x)|
\leq c_{26} \cdot \frac{a_n^p}{M_n^p}
\leq c_{26} \cdot (\log n) \cdot \frac{1}{M_n^p}.
\]
Plugging in the values for $R_n$ and $M_n$ we get the assertion.
\quad \hfill $\Box$

\subsection{Proof of Theorem \ref{th2}}
W.l.o.g. we assume $supp(X) \subseteq [-A_n,A_n]^d$.

Define the estimate $\bar{m}_n$ exactly like $m_n$ except that
for given directions $\bb_l$ $(l=1,\dots,r)$
we define the neural network estimate $\tilde{m}_n(x)$ by
\[
\tilde{m}_n(x)
=
\sum_{l=1,\dots,r}
\sum_{k=1}^{M_n+1}
\sum_{\substack{j_1, \dots, j_d \in \{0, \dots, N\}\\j_1 + \dots+j_d \leq N}}
a_{i_k,j_1,\dots,j_d,\bb_l}
\cdot
f_{net,j_1,\dots,j_d,i_k,\bb_l}(x),
\]
where the coefficients $a_{k ,j_1,\dots,j_d,\bb_l}$ are chosen
from the set
\[
\left\{
(a_{k ,j_1,\dots,j_d,\bb_l})_{k,j_1, \dots, j_d,l} \, : \,
\sum_{k,j_1,\dots,j_d,l}
a_{k ,j_1,\dots,j_d,\bb_l}^2
\leq
c_{27} \cdot n^2
\right\}
\]
by minimizing
\[
\frac{1}{n} \sum_{i=1}^n |Y_i-\tilde{m}_n(X_i)|^2
+
\frac{c_3}{n} \cdot
\sum_{l=1}^r
\sum_{k=0}^{K}
\sum_{\substack{j_1, \dots, j_d \in \{0, \dots, N\}\\j_1 + \dots+j_d \leq N}}
a_{k,j_1,\dots,j_d,\bb_l}^2
\]
for some constant $c_3>0$. Then $\bar{m}_n$ satisfies
\[
\bar{m}_n \in \left\{ T_{\beta_n} f \, : \, f \in \F^{(J_n)} \right\}, 
\]
where $\F^{(J_n)}$ (with $J_n=r \cdot (M_n+1) \cdot \left( {N+d \atop
    d}\right)$)
 is the function space defined in Lemma \ref{le8}.
On the event
\[
B_n = \{ |Y_i| \leq \sqrt{n} \, : \, i=1, \dots, n\}
\]
we know by (\ref{se3eq9b}) that we have $m_n=\bar{m}_n$
(provided $c_{27} \geq 1/c_3$). Hence
\[
\int |m_n(x)-m(x)|^2 \PROB_X(dx)
\leq
\int |\bar{m}_n(x)-m(x)|^2 \PROB_X(dx) + 4 \beta_n^2 \cdot 1_{B_n^c}.
\]
By Markov inequality we know
\[
\PROB\{B_n^c\} \leq n \cdot \PROB\{|Y|> \sqrt{n}\}
\leq
\frac{n \cdot \EXP\{ e^{c_3 \cdot Y^2} \} }{\exp(c_3 \cdot n)},
\]
therefore (\ref{subgaus}) implies that it suffices to show the assertion
under the additional assumption
\begin{equation}
\label{pth2eq1}
\tilde{m}_n(\cdot, (X_1,Y_1), \dots, (X_n,Y_n))
\in \F^{(J_n)}.
\end{equation}
By Lemma \ref{le2b} we know that for each
$i \in \{1, \dots, I_n\}$ there exist coefficients
$a_{k,j_1,\dots,j_d,l}^{(i)} \in [-c_{28} \cdot A_n^p ,c_{28} \cdot A_n^p]$, which depend
on $\ba_l$ and on $\bb_l^{(i)}$, but which are independent of
$(X_1,Y_1)$, \dots, $(X_n,Y_n)$, such that
we have for all $x \in [-A_n,A_n]^d$
\begin{eqnarray}
	&&
	\Bigg|
	\sum_{l=1}^r \sum_{j=0}^q \frac{
		g_l^{(j)}((\bb_l^{(i)})^Tx)
	}{j!}
	\cdot
	((\ba_l-\bb_l^{(i)})^T x)^j \nonumber \\
	&&
	\hspace*{2cm}
	-
	\sum_{l=0}^r
	\sum_{k=1}^{M_n+1}
	\sum_{\substack{j_1, \dots, j_d \in \{0, \dots, N\}\\j_1 + \dots+j_d \leq N}}
	a_{i_k,j_1,\dots,j_d,l}^{(i)} \cdot (x^{(1)})^{j_1} \cdots (x^{(d)})^{j_d}
	\nonumber \\
	&&
	\hspace*{4cm}
	\cdot \left(
	1 -  \frac{M_n}{2 \cdot \sqrt{d} \cdot A_n} \cdot |(\bb_l^{(i)})^Tx - u_{i_k}|
	\right)_+
	\Bigg|
	\nonumber \\
	&&
	\leq
	r \cdot 2^{p} \cdot (p+1) \cdot C  \cdot A_n^{2p} \cdot
	\left(
	\max \left\{
	\frac{1}{M_n}, \max_{l=1, \dots,r }\|\ba_l-\bb_l^{(i)}\|_\infty
	\right\}
	\right)^p.
        \label{pth2eq2}
\end{eqnarray}  
From the definition of the estimate 
we get
\begin{eqnarray*}
	&&
	\frac{1}{n} \sum_{i=1}^n |Y_i-\tilde{m}_n(X_i)|^2
	\\
	&&
	\leq
	\min_{t=1, \dots, I_n}
	\Bigg\{
	\frac{1}{n} \sum_{i=1}^n |Y_i
	-
	\sum_{l=1,\dots,r}
	\sum_{k=1}^{M_n+1}
	\sum_{\substack{j_1, \dots, j_d \in \{0, \dots, N\}\\j_1 + \dots+j_d \leq N}}
	a_{i_k,j_1,\dots,j_d,l}^{(t)}
	\cdot
	f_{net,j_1,\dots,j_d,i_k,\bb_l^{(t)}}(X_i)
	|^2
	\\
	&&
	\hspace*{3cm}
	+
	\frac{c_3}{n} \cdot
	\sum_{l=1}^r
	\sum_{k=1}^{M_n+1}
	\sum_{\substack{j_1, \dots, j_d \in \{0, \dots, N\}\\j_1 + \dots+j_d \leq N}}
	(a_{i_k,j_1,\dots,j_d,l}^{(t)})^2
	\Bigg\}
	\\
	&&
	\leq
	\min_{t=1, \dots, I_n}
	\Bigg\{
	\frac{1}{n} \sum_{i=1}^n |Y_i
	-
	\sum_{l=1,\dots,r}
	\sum_{k=1}^{M_n+1}
	\sum_{\substack{j_1, \dots, j_d \in \{0, \dots, N\}\\j_1 + \dots+j_d \leq N}}
	a_{i_k,j_1,\dots,j_d,l}^{(t)}
	\cdot
	f_{net,j_1,\dots,j_d,i_k,\bb_l^{(t)}}(X_i)
	|^2
	\\
	&&
	\quad
	+
	c_{29} \cdot A_n^{2p} \cdot r \cdot \left( {N+d \atop d} \right) \cdot \frac{M_n}{n} \Bigg\}.
\end{eqnarray*}
Hence,
application of Lemma \ref{le6}
and Lemma \ref{le8}
(together with (\ref{pth2eq1}))
yields
\begin{eqnarray*}
	&&
	\EXP \int |m_n(x)-m(x)|^2 \PROB_X (dx)
	\\
	&&
	\leq
	c_{30} \cdot \frac{(\log n)^3 \cdot M_n}{n}
	\\
	&&
	\quad
	+
	2 \cdot
	\EXP \Bigg(
	\min_{t=1,\dots,I_n}
	\int
	\Big|
	\sum_{l=1,\dots,r}
	\sum_{k=1}^{M_n+1}
	\sum_{\substack{j_1, \dots, j_d \in \{0, \dots, N\}\\j_1 + \dots+j_d \leq N}}
	a_{i_k,j_1,\dots,j_d,l}^{(t)}
	\cdot
	f_{net,j_1,\dots,j_d,i_k,\bb_l^{(t)}}(x)
	\\
	&&
	\hspace*{6cm}
	-m(x)
	\Big|^2 \PROB_{X}(dx)\Bigg)
	+
	c_{31} \cdot (\log n) \cdot n^{-\frac{2p}{2p+1}}.
\end{eqnarray*}
Because of $(a+b+c)^2 \leq 3 a^2 + 3 b^2 + 3 c^2$ $(a,b,c \in \R)$
we have 
\begin{eqnarray*}
	&&
	\int
	|
	\sum_{l=1,\dots,r}
	\sum_{k=1}^{M_n}
	\sum_{\substack{j_1, \dots, j_d \in \{0, \dots, N\}\\j_1 + \dots+j_d \leq N}}
	a_{i_k,j_1,\dots,j_d,l}^{(t)}
	\cdot
	f_{net,j_1,\dots,j_d,i_k,\bb_l^{(t)}}(x)
	-m(x)
	|^2 \PROB_X(dx)
	\\
	&&
	\leq
	3 \cdot
	\int
	\Bigg|
	\sum_{l=1,\dots,r}
	\sum_{k=1}^{M_n}
	\sum_{\substack{j_1, \dots, j_d \in \{0, \dots, N\}\\j_1 + \dots+j_d \leq N}}
	a_{i_k,j_1,\dots,j_d,l}^{(t)}
	\cdot
	f_{net,j_1,\dots,j_d,i_k,\bb_l^{(t)}}(x)
	\\
	&&
	\hspace*{1cm}
	-
	\sum_{l=1}^r \sum_{j=0}^q
	\sum_{\substack{j_1, \dots, j_d \in \{0, \dots, N\}\\j_1 + \dots+j_d \leq N}}
	a_{i_k,j_1,\dots,j_d,l}^{(t)} \cdot (x^{(1)})^{j_1} \cdots (x^{(d)})^{j_d}
	\\
	&&
	\hspace*{4cm}
	\cdot \left(
	1 -  \frac{M_n}{2 \cdot \sqrt{d} \cdot A_n} \cdot |(\bb_l^{(t)})^Tx - u_{i_k}|
	\right)_+
	\Bigg|^2 \PROB_X(dx)
	\\
	&&
	\quad
	+
	3 \cdot
	\int
	|
	\sum_{l=1}^r \sum_{j=0}^q
	\sum_{\substack{j_1, \dots, j_d \in \{0, \dots, N\}\\j_1 + \dots+j_d \leq N}}
	a_{i_k,j_1,\dots,j_d,l}^{(t)} \cdot (x^{(1)})^{j_1} \cdots (x^{(d)})^{j_d}
	\\
	&&
	\hspace*{4cm}
	\cdot \left(
	1 -  \frac{M_n}{2 \cdot \sqrt{d} \cdot A_n} \cdot |(\bb_l^{(t)})^Tx - u_{i_k}|
	\right)_+
	\\
	&&
	\hspace*{4cm}
	-
	\sum_{l=1}^r \sum_{j=0}^q \frac{
		g_l^{(j)}((\bb_l^{(t)})^Tx)
	}{j!}
	\cdot
	((\ba_l-\bb_l^{(t)})^T x)^j
	|^2 \PROB_X(dx)
	\\
	&&
	\quad
	+
	3 \cdot
	\int
	|
	\sum_{l=1}^r \sum_{j=0}^q \frac{
		g_l^{(j)}(\bb_l^Tx)
	}{j!}
	\cdot
	((\ba_l-\bb_l^{(t)})^T x)^j
	-
	m(x)
	|^2 \PROB_X (dx).
\end{eqnarray*}
Application of Lemma \ref{le5}
implies for all $x \in [-A_n,A_n]^d$
\begin{eqnarray*}
	&&
	\Bigg|
	\sum_{l=1,\dots,r}
	\sum_{k=0}^{M_n}
	\sum_{\substack{j_1, \dots, j_d \in \{0, \dots, N\}\\j_1 + \dots+j_d \leq N}}
	a_{k,j_1,\dots,j_d,l}^{(t)}
	\cdot
	f_{net,k,j_1,\dots,j_d,\bb_l^{(t)}}(x)
	\\
	&&
	\hspace*{1cm}
	-
	\sum_{l=1}^r \sum_{j=0}^q
	\sum_{\substack{j_1, \dots, j_d \in \{0, \dots, N\}\\j_1 + \dots+j_d \leq N}}
	a_{k,j_1,\dots,j_d,l}^{(t)} \cdot (x^{(1)})^{j_1} \cdots (x^{(d)})^{j_d}
	\\
	&&
	\hspace*{4cm}
	\cdot \left(
	1 -  \frac{K}{2 \cdot \sqrt{d} \cdot A} \cdot |(\bb_l^{(t)})^Tx - u_k|
	\right)_+
	\Bigg|^2
	\\
	&&
	\leq r^2 \cdot (M_n+1)^2 \cdot (N+d)^{2d} 
\cdot c_{28}^2 \cdot A_n^{2p}
\cdot c_{12}^2
        \cdot 3^{6 \cdot 3^s} \cdot A_n^{6 \cdot 2^s} 
	\cdot M_n^6 \cdot \frac{1}{R_n^2}
	\leq c_{32} \cdot \frac{(\log n)^2}{n}. 
\end{eqnarray*}
By Lemma \ref{le1b} we have for all $x \in [-A_n,A_n]^d$
\begin{eqnarray*}
	&&
	|
	\sum_{l=1}^r \sum_{j=0}^q \frac{
		g_l^{(j)}((\bb_l^{(t)})^Tx)
	}{j!}
	\cdot
	((\ba_l-\bb_l)^T x)^j
	-
	m(x)
	|^2
	\leq
	c_{33} \cdot A_n^{2p} \cdot \|\ba_l-\bb_l^{(t)}\|_\infty^{2p}.
\end{eqnarray*}
Using this together with (\ref{pth2eq2}) we see that
it remains to show 
\[
\EXP
\left\{
\min_{i=1,\dots,I_n}
\max_{s=1, \dots,r}
\|\bb_s^{(i)} - \ba_s\|_\infty^{2p}
\right\}
\leq
c_{34}
\cdot
(\log n)^2 \cdot n^{-\frac{2p}{2p+1}}
.
\]
By the random choice of the 
$\bb_l^{(i)}$
we know for any $t \in (0,1]$
\begin{eqnarray*}
	\PROB\left\{
	\min_{i=1,\dots,I_n}
	\max_{l=1, \dots,r}
	\|\bb_l^{(i)} - \ba_l\|_\infty
	>t
	\right\}
	&=&
	\prod_{i=1}^{I_n} (1-\PROB\left\{
	\max_{l=1, \dots,r}
	\|\bb_l^{(i)} - \ba_l\|_\infty
	\leq t
	\right\})
	\\
	&\leq&
	\left(
	1-\left(\frac{t}{2}\right)^{r \cdot d}
	%1- t^{r \cdot d}
	\right)^{I_n}
\end{eqnarray*}
from which we conclude
\begin{eqnarray*}
	&&
	\EXP
	\left\{
	\min_{i=1,\dots,I_n}
	\max_{l=1, \dots,r}
	\|\bb_l^{(i)} - \ba_l \|_\infty^{2p}
	\right\}
	\\
	&&
	\leq
	\left(
	\frac{\log n}{n}
	\right)^{\frac{2p}{2p+1}}
	+
	2^{2p} \cdot
	\PROB\left\{
	\min_{i=1,\dots,I_n}
	\max_{l=1, \dots,r}
	\|\bb_l^{(i)} - \ba_l\|_\infty
	>
	\left(
	\frac{\log n}{n}
	\right)^{\frac{1}{2p+1}}
	\right\}
	\\
	&&
	\leq
	\left(
	\frac{\log n}{n}
	\right)^{\frac{2p}{2p+1}}
	+
	2^{2p} \cdot
	\left(
	1
	-
\cdot \frac{1}{2^{r \cdot d}}
	\left(
	\frac{\log n }{n}
	\right)^{\frac{r \cdot d}{2p+1}}
	\right)^{I_n}
\\
&&
\leq
	c_{35}
	\cdot
\exp
	\left(
- I_n \cdot
\frac{1}{2^{r \cdot d}}
\cdot
	\left(
	\frac{\log n }{n}
	\right)^{\frac{r \cdot d}{2p+1}}
\right)
	\\
&&
=
	c_{35}
	\cdot
\exp
	\left(
- \frac{c_9}{2^{r \cdot d}} \cdot (\log n)^2
\right) 
	\\
	&&
	\leq
	c_{35}
	\cdot
	\left(
	\frac{\log n}{n}
	\right)^{\frac{2p}{2p+1}}
\end{eqnarray*}
where the last inequality follows from
\[
I_n \geq
c_9 \cdot
(\log n)^2 \cdot
\left(
\frac{n}{\log n}
\right)^{\frac{r \cdot d}{2p+1}}.
\]
Putting together the above results we get the assertion.
\hfill $\Box$

\newpage
\begin{appendix}
  
\section{Supplementary material}

\subsection{Computation of the linear neural network estimate}
The estimate in Subsection \ref{se2sub2}  is given by
\begin{equation}
\label{se2eq8}
\tilde{m}_n(x)
=
\sum_{j=1}^J
a_j \cdot B_j(x)
\end{equation}
where ${\mathbf a}=(a_j)_{j=1,\dots,J} \in \R^J$ minimizes
\begin{eqnarray*}
&&
\frac{1}{n} ({\mathbf Y}- {\mathbf B a})^T ({\mathbf Y} - {\mathbf B
  a})
+
\frac{c_3}{n} \cdot{\mathbf a}^T {\mathbf a}
\\
&&
=
\frac{1}{n} ({\mathbf Y}^T {\mathbf Y}
-
2 {\mathbf Y}^T  {\mathbf B
  a})
+
{\mathbf  a}^T
\left(
\frac{1}{n}
 {\mathbf B}^T  {\mathbf B} +
\frac{c_3}{n} \cdot {\mathbf 1}
\right)
{\mathbf  a}.
\end{eqnarray*}
Since the matrix
\[
{\mathbf A}
=
\frac{1}{n}
 {\mathbf B}^T  {\mathbf B} +
\frac{c_3}{n} \cdot {\mathbf 1}
\]
is positive definite, its inverse matrix ${\mathbf A}^{-1}$ exists
and it is easy to see that we have
\begin{eqnarray*}
&&
\frac{1}{n} \cdot ( {\mathbf Y}^T {\mathbf Y}
-
2 {\mathbf Y}^T  {\mathbf B
  a})
+
{\mathbf  a}^T
\left(
\frac{1}{n}
 {\mathbf B}^T {\mathbf B} +
\frac{c_3}{n} \cdot {\mathbf 1}
\right)
{\mathbf  a}
\\
&&
=
( {\mathbf a} -\frac{1}{n} \cdot  {\mathbf A}^{-1}  {\mathbf B}^T {\mathbf Y})^T
{\mathbf A}
( {\mathbf a} - \frac{1}{n} \cdot {\mathbf A}^{-1}  {\mathbf B}^T {\mathbf Y})
+ \frac{1}{n} {\mathbf Y}^T {\mathbf Y}
-
\frac{1}{n^2} \cdot
{\mathbf Y}^T {\mathbf B} {\mathbf A}^{-1} {\mathbf B}^T {\mathbf Y}.
\end{eqnarray*}
The last expression is minimal for $ {\mathbf a} = \frac{1}{n} \cdot {\mathbf A}^{-1}
{\mathbf B}^T {\mathbf Y}$,
which proves that the vector of coefficients of our estimate (\ref{se2eq8})
is the unique solution of the linear equation system (\ref{se2eq9}).

\subsection{Proof of Lemma \ref{le5}}
    
Define $g_1^{(0)}$ by backward recursion:
\[
g_k^{(s)}(x)=x^{(l)}-y^{(l)}
\]
for $j_1 + j_2 + \dots + j_{l-1}+1 \leq k \leq j_1 + j_2 + \dots +
j_l$
and $l=1,\dots,d$,
\[
g^{(s)}_{j_1+j_2+\dots+j_d+k}(x) =
\left(
1 - \frac{M}{2a} \cdot |x^{(k)}-y^{(k)}|
\right)_+
\]
for $k=1, \dots, d$, and
\[
g^{(s)}_k(x)=1
\]
for $k=j_1+j_2+\dots+j_d+d+1, j_1+j_2+\dots+j_d+d+2, \dots,2^s$,
and
\[
g_k^{(l)}(x)
=
g_{2k-1}^{(l+1)}(x) \cdot g_{2k}^{(l+1)}(x)
\]
for $k \in \{1, 2, \dots, 2^l\}$ and $l \in \{0,\dots,s-1\}$.

Then we have for any $l \in \{0,\dots,s\}$, $k \in \{1,\dots,2^l\}$
and
$x \in [-a,a]$
\[
|g_k^{(l)}(x)| \leq (2a)^{2^{s-l}}.
\]
By Lemma \ref{le2} the network $f_{mult}$ satisfies for any
$l \in \{0, \dots, s\}$ and \linebreak
$x,y \in [-3^{3^{s-l}} \cdot a^{2^{s-l}}, 3^{3^{s-l}} \cdot a^{2^{s-l}}]$
\[
|f_{mult}(x,y)-x \cdot y| \leq \frac{20 \cdot \| \sigma^{\prime \prime \prime}\|_\infty}{
3 \cdot |\sigma^{\prime\prime}(t_{\sigma})|}
\cdot 3^{3 \cdot 3^{s-l}} \cdot a^{3 \cdot 2^{s-l}}
\cdot \frac{1}{R}.
\]
Furthermore we have by Lemma \ref{le1} and Lemma \ref{le4}
for any $x \in [-3a,3a]$
\begin{equation}
\label{ple5eq1}
| f_{id}(x)-x|
\leq
\frac{9 \cdot
\| \sigma^{\prime \prime}\|_{\infty} \cdot a^2
}{
2 \cdot |\sigma^\prime(t_{\sigma,id})|
}
\cdot
\frac{1}{R}
\end{equation}
and for any $x \in [-a,a]^d$
\begin{eqnarray}
\label{ple5eq2} \nonumber
&&
\left| f_{hat,y}(x)-
\left(
1 - \frac{M}{2a}  \cdot |x-y|
\right)_+
\right|
\\
&&
\leq
1792 \cdot \frac{
\max \left\{
\| \sigma^{\prime \prime}\|_{\infty},
\| \sigma^{\prime \prime \prime}\|_{\infty},1
\right\}
}{
\min \left\{
2 \cdot |\sigma^\prime (t_{\sigma.id})|,
|\sigma^{\prime \prime} (t_{\sigma})|, 1
\right\}
} \cdot M^3  \cdot \frac{1}{R}.
\end{eqnarray}
From this and (\ref{le5eq1}) we can recursively conclude
\[
|f^{(l)}_k(x)|
\leq
3^{3^{s-l}}
\cdot
a^{2^{s-l}}
\]
for $k \in \{1, \dots, 2^l\}$ and $l \in  \{0,\dots,s\}$.

In order to prove the assertion of Lemma \ref{le5} we show in the
sequel
\[
| f_k^{(l)}(x)-g_k^{(l)}(x)|
\leq
c_{36} \cdot
3^{3 \cdot 3^{s-l}} \cdot a^{3 \cdot 2^{s-l}} \cdot M^3 \cdot \frac{1}{R}
\]
for $k \in \{1, \dots, 2^l\}$ and $l \in  \{0,\dots,s\}$,
where
\[
c_{36}
=
\max \left\{
\frac{20 \cdot \| \sigma^{\prime \prime \prime}\|_\infty}{
  3 \cdot |\sigma^{\prime\prime}(t_{\sigma})|},
\frac{9 \cdot
\| \sigma^{\prime \prime}\|_{\infty} 
}{
|\sigma^\prime(t_{\sigma,id})|
}
,
1792 \cdot \frac{
\max \left\{
\| \sigma^{\prime \prime}\|_{\infty},
\| \sigma^{\prime \prime \prime}\|_{\infty},1
\right\}
}{
\min \left\{
2 \cdot |\sigma^\prime (t_{\sigma.id})|,
|\sigma^{\prime \prime} (t_{\sigma})|, 1
\right\}
} 
\right\}
\]

For $s=l$ this is a consequence of (\ref{ple5eq1}),
and (\ref{ple5eq2}). For $l \in \{0,1, \dots, s-1\}$
we can conclude via induction
\begin{eqnarray*}
&&
| f_{k}^{(l)}(x)-g_k^{(l)}(x)|
\\
&&
\leq
| f_{mult}(f_{2k-1}^{(l+1)}(x),f_{2k}^{(l+1)}(x))
-
f_{2k-1}^{(l+1)}(x) \cdot f_{2k}^{(l+1)}(x)
|
\\
&&
\quad
+
|f_{2k-1}^{(l+1)}(x) \cdot f_{2k}^{(l+1)}(x)
-
g_{2k-1}^{(l+1)}(x) \cdot f_{2k}^{(l+1)}(x)|
\\
&&
\quad
+
|g_{2k-1}^{(l+1)}(x) \cdot f_{2k}^{(l+1)}(x)
-
g_{2k-1}^{(l+1)}(x) \cdot g_{2k}^{(l+1)}(x)|
\\
&&
\leq
c_{36}
\cdot 3^{3 \cdot 3^{s-l-1}} \cdot a^{3
  \cdot 2^{s-l-1}}
\cdot \frac{1}{R}
+
3^{3^{s-l-1}}
\cdot
a^{2^{s-l-1}}
\cdot
2
\cdot
c_{36} \cdot
3^{3 \cdot 3^{s-l-1}} \cdot a^{3 \cdot 2^{s-l-1}} \cdot M^3 \cdot \frac{1}{R}
\\
&&
\leq
c_{36}
 \cdot
 \left(
3^{3^{s-l}} + 2 \cdot  3^{4 \cdot 3^{s-l-1}} 
\right)
\cdot a^{3 \cdot 2^{s-l}} \cdot M^3 \cdot \frac{1}{R}
\\
&&
\leq
c_{36}
 \cdot
3^{3 \cdot 3^{s-l}} \cdot a^{3 \cdot 2^{s-l}} \cdot M^3 \cdot \frac{1}{R}
.
\end{eqnarray*}
\quad \hfill $\Box$

\subsection{Lemma \ref{le5} in the context of projection pursuit}

\begin{lemma}
	\label{le5projpurs} 
	Let $M \in \N$ and
	let $\sigma: \mathbb{R} \to [0,1]$ be 2-admissible  according to
	Definition \ref{se5de1}. Let $A \geq 1$,
	$\bb \in \Rd$ with $\|\bb\| \leq 1$ and
	\begin{eqnarray}
	\label{le7eq1}
	R
	& \geq&
	\max\Bigg\{
	\frac{
		\| \sigma^{\prime \prime}\|_{\infty}
		\cdot  (M+1)
	}{
		2 \cdot  | \sigma^\prime (t_{\sigma,id})|
	}
	,
	\frac{
		9 \cdot
		\| \sigma^{\prime \prime}\|_{\infty}\cdot A
	}{
		| \sigma^\prime (t_{\sigma,id})|
	}
	,\\
	&&
	\frac{20 \cdot \| \sigma^{\prime \prime \prime}\|_\infty}{
		3 \cdot |\sigma^{\prime\prime}(t_{\sigma})|}
	\cdot 3^{3 \cdot 3^s} \cdot A^{3 \cdot 2^s},
	1792 \cdot \frac{
		\max \left\{
		\| \sigma^{\prime \prime}\|_{\infty},
		\| \sigma^{\prime \prime \prime}\|_{\infty},1
		\right\}
	}{
		\min \left\{
		2 \cdot |\sigma^\prime (t_{\sigma,id})|,
		|\sigma^{\prime \prime} (t_{\sigma})|, 1
		\right\}
	} \cdot d^{3/2} \cdot M^3 
	\Bigg\}
	\nonumber
	\end{eqnarray}
	and let $y \in [-A,A]$. Let $N \in \N$ and let
	$j_1, \dots, j_d \in \N_0$ such that $j_1+\dots+j_d \leq N$, and set
	$s=\lceil \log_2(N+1) \rceil$.
	Let $f_{id}$, $f_{mult}$ and $\bar{f}_{hat,z}$ (for $z \in \R$) be the
	neural
	networks defined in Subsection \ref{se3sub2}. (So in particular
	$\bar{f}_{hat,z}$ is the neural network from Lemma \ref{le4}
	with $y=z$ and $a=\sqrt{d} \cdot A$.)
	Define the network $f_{net,j_1,\dots,j_d,y}$ by
	\[
	f_{net,j_1,\dots,j_d,y}(x)=f_1^{(0)}(x),
	\]
	where $f_1^{(0)}$ is defined by backward recursion as follows:
	\[
	f_k^{(l)}(x)
	=
	f_{mult}
	\left(
	f_{2k-1}^{(l+1)}(x),f_{2k}^{(l+1)}(x)
	\right)
	\]
	for $k \in \{1, 2, \dots, 2^l\}$ and $l \in \{0,\dots,s-1\}$,
	and
	\[
	f_k^{(s)}(x)=f_{id}(f_{id}(x^{(l)}))
	\]
	for $j_1 + j_2 + \dots + j_{l-1}+1 \leq k \leq j_1 + j_2 + \dots +
	j_l$
	and $l=1,\dots,d$,
	\[
	f_{j_1+j_2+\dots+j_d+1}^{(s)}(x) = \bar{f}_{hat,y}(\bb^T x),
	\]
	and
	\[
	f^{(s)}_k(x)=1
	\]
	for $k=j_1+j_2+\dots+j_d+2, j_1+j_2+\dots+j_d+3, \dots,2^s$.
	Then we have for any $x \in [-A,A]^d$:
	\begin{eqnarray*}
		&&
		\left|
		f_{net,y}(x)
		-
		(x^{(1)})^{j_1} \cdots (x^{(d)})^{j_d}
		\cdot (1 - \frac{M}{2 \cdot \sqrt{d} \cdot A}
		\cdot |\bb^T x - y|)_+
		\right|
		\\
		&&
		\leq
		c_{37} \cdot 3^{3 \cdot 3^s} \cdot A^{3 \cdot 2^s} 
		\cdot M^3 \cdot \frac{1}{R}.
	\end{eqnarray*}
\end{lemma}

\noindent
{\bf Proof.}
Define $g_1^{(0)}$ by backward recursion:
\[
g_k^{(s)}(x)=x^{(l)}
\]
for $j_1 + j_2 + \dots + j_{l-1}+1 \leq k \leq j_1 + j_2 + \dots +
j_l$
and $l=1,\dots,d$,
\[
g^{(s)}_{j_1+j_2+\dots+j_d+1}(x) =
\left(
1 - \frac{M}{2 \cdot \sqrt{d} \cdot A} \cdot |\bb^T x-y|
\right)_+,
\]
and
\[
g^{(s)}_k(x)=1
\]
for $k=j_1+j_2+\dots+j_d+2, j_1+j_2+\dots+j_d+3, \dots,2^s$,
and
\[
g_k^{(l)}(x)
=
g_{2k-1}^{(l+1)}(x) \cdot g_{2k}^{(l+1)}(x)
\]
for $k \in \{1, 2, \dots, 2^l\}$ and $l \in \{0,\dots,s-1\}$.

Then we have for any $x \in [-A,A]^d$
\[
|g_k^{(l)}(x)| \leq A^{2^{s-l}}.
\]
By Lemma \ref{le2} the network $f_{mult}$ satisfies for any
$l \in \{0, \dots, s\}$ and \linebreak
$x,y \in [-3^{3^{s-l}} \cdot A^{2^{s-l}}, 3^{3^{s-l}} \cdot A^{2^{s-l}}]$
\[
|f_{mult}(x,y)-x \cdot y| \leq \frac{20 \cdot \| \sigma^{\prime \prime \prime}\|_\infty}{
	3 \cdot |\sigma^{\prime\prime}(t_{\sigma})|}
\cdot 3^{3 \cdot 3^{s-l}} \cdot A^{3 \cdot 2^{s-l}}
\cdot \frac{1}{R}.
\]
Furthermore we have by Lemma \ref{le1} and Lemma \ref{le4}
for any $x \in [-3A,3A]$
\begin{equation}
\label{pl7eq1}
| f_{id}(x)-x|
\leq
\frac{9 \cdot
	\| \sigma^{\prime \prime}\|_{\infty} \cdot A^2
}{
	2 \cdot |\sigma^\prime(t_{\sigma,id})|
}
\cdot
\frac{1}{R}
\end{equation}
and for any $x \in [-A,A]^d$
\begin{eqnarray}
\label{pl7eq2} \nonumber
&&
\left| \bar{f}_{hat,y}(x)-
\left(
1 - \frac{M}{2 \cdot \sqrt{d} \cdot A}  \cdot |\bb^T x-y|
\right)_+
\right|
\\
&&
\leq
1792 \cdot \frac{
	\max \left\{
	\| \sigma^{\prime \prime}\|_{\infty},
	\| \sigma^{\prime \prime \prime}\|_{\infty},1
	\right\}
}{
	\min \left\{
	2 \cdot |\sigma^\prime (t_{\sigma.id})|,
	|\sigma^{\prime \prime} (t_{\sigma})|, 1
	\right\}
}  \cdot M^3  \cdot \frac{1}{R}.
\end{eqnarray}
From this and (\ref{le7eq1}) we can recursively conclude
\[
|f^{(l)}_k(x)|
\leq
3^{3^{s-l}}
\cdot
A^{2^{s-l}}
\]
for $k \in \{1, \dots, 2^l\}$ and $l \in  \{0,\dots,s\}$.

In order to prove the assertion of Lemma \ref{le5} we show in the
sequel
\[
| f_k^{(l)}(x)-g_k^{(l)}(x)|
\leq
c_{37} \cdot
3^{3 \cdot 3^{s-l}} \cdot A^{3 \cdot 2^{s-l}} 
\cdot M^3 \cdot \frac{1}{R}
\]
for $k \in \{1, \dots, 2^l\}$ and $l \in  \{0,\dots,s\}$,
where
\[
c_{37}
=
\max \left\{
\frac{20 \cdot \| \sigma^{\prime \prime \prime}\|_\infty}{
	3 \cdot |\sigma^{\prime\prime}(t_{\sigma})|},
\frac{9 \cdot
	\| \sigma^{\prime \prime}\|_{\infty} 
}{
	|\sigma^\prime(t_{\sigma,id})|
}
,
1792 \cdot \frac{
	\max \left\{
	\| \sigma^{\prime \prime}\|_{\infty},
	\| \sigma^{\prime \prime \prime}\|_{\infty},1
	\right\}
}{
	\min \left\{
	2 \cdot |\sigma^\prime (t_{\sigma.id})|,
	|\sigma^{\prime \prime} (t_{\sigma})|, 1
	\right\}
}
\right\}.
\]

For $s=l$ this is a consequence of (\ref{pl7eq1}) and
(\ref{pl7eq2}). For $l \in \{0,1, \dots, s-1\}$
we can conclude via induction
\begin{eqnarray*}
	&&
	| f_{k}^{(l)}(x)-g_k^{(l)}(x)|
	\\
	&&
	\leq
	| f_{mult}(f_{2k-1}^{(l+1)}(x),f_{2k}^{(l+1)}(x))
	-
	f_{2k-1}^{(l+1)}(x) \cdot f_{2k}^{(l+1)}(x)
	|
	\\
	&&
	\quad
	+
	|f_{2k-1}^{(l+1)}(x) \cdot f_{2k}^{(l+1)}(x)
	-
	g_{2k-1}^{(l+1)}(x) \cdot f_{2k}^{(l+1)}(x)|
	\\
	&&
	\quad
	+
	|g_{2k-1}^{(l+1)}(x) \cdot f_{2k}^{(l+1)}(x)
	-
	g_{2k-1}^{(l+1)}(x) \cdot g_{2k}^{(l+1)}(x)|
	\\
	&&
	\leq
	c_{37}
	\cdot 3^{3 \cdot 3^{s-l-1}} \cdot A^{3
		\cdot 2^{s-l-1}}
	\cdot \frac{1}{R}
	+
	3^{3^{s-l-1}}
	\cdot
	A^{2^{s-l-1}}
	\cdot
	2
	\cdot
	c_{37} \cdot
	3^{3 \cdot 3^{s-l-1}} \cdot A^{3 \cdot 2^{s-l-1}}  \cdot  M^3
	\cdot \frac{1}{R}
	\\
	&&
	\leq
	c_{37}
	\cdot
	\left(
	3^{3^{s-l}} + 2 \cdot  3^{4 \cdot 3^{s-l-1}} 
	\right)
	\cdot A^{3 \cdot 2^{s-l}}  \cdot M^3 \cdot \frac{1}{R}
	\\
	&&
	\leq
	c_{37}
	\cdot
	3^{3 \cdot 3^{s-l}} \cdot A^{3 \cdot 2^{s-l}}  \cdot M^3 \cdot \frac{1}{R}
	.
\end{eqnarray*}
\quad \hfill $\Box$

\subsection{Proof of Lemma \ref{le6}}
In the proof we use the following error decomposition:
\beq
        \int |m_n(x) - m(x)|^2 \PROB_X (dx)&& \qquad \qquad\qquad \qquad\qquad \qquad\qquad \qquad\qquad \qquad
\eeq
\vspace{-8mm}
\beq
        &=&\Big[ \textbf{E}\Big\{|m_n(X) -Y|^2|\mathcal{D}_n\Big\} - \textbf{E} \Big\{|m(X)-Y|^2 \Big\} \\
  &&\qquad -
\Big(
  \textbf{E} \Big\{|m_n(X) - T_{\beta_n} Y|^2|\mathcal{D}_n \Big\} - \textbf{E} \Big\{|m_{\beta_n}(X)-
        T_{\beta_n} Y|^2 \Big\} \Big) \Big]\\
        &&+\Bigg[ \textbf{E} \Big\{|m_n(X) - T_{\beta_n} Y|^2|\mathcal{D}_n \Big\} - \textbf{E} \Big\{|m_{\beta_n}(X)
        - T_{\beta_n} Y|^2 \Big\}\\
        &&\qquad  - 2 \cdot \frac{1}{n} \sum_{i=1}^n \Big(|m_n(X_i) - T_{\beta_n} Y_i|^2 - |m_{\beta_n} (X_i)
  - T_{\beta_n} Y_i|^2 \Big) \Bigg]\\
        && + \left[ 2 \cdot \frac{1}{n} \sum_{i=1}^n |m_n(X_i) - T_{\beta_n} Y_i|^2
                                                - 2\cdot \frac{1}{n} \sum_{i=1}^n |m_{\beta_n} (X_i) - T_{\beta_n} Y_i|^2 \right. \\
        && \qquad \left. -\left( 2\cdot \frac{1}{n} \sum_{i=1}^n |m_n(X_i)-Y_i|^2 - 2\cdot \frac{1}{n} \sum_{i=1}^n|m(X_i)-Y_i|^2 \right)\right]\\
        && + \left[2\left( \frac{1}{n} \sum_{i=1}^n |m_n(X_i)-Y_i|^2 -\frac{1}{n} \sum_{i=1}^n|m(X_i)-Y_i|^2\right)\right]\\
        &=& \sum_{i=1}^4 T_{i,n},
\eeq
where $T_{\beta_n} Y$ is the truncated version of $Y$ and $m_{\beta_n}$ is the regression function of $T_{\beta_n} Y$, i.e.,
\[ m_{\beta_n}(x) = \textbf {E} \Big\{ T_{\beta_n} Y|X=x\Big\}.\]
We start with bounding $T_{1,n}$. By using $a^2-b^2 =(a-b)(a+b)$ we get
\beq
        T_{1,n}&=& \textbf{E} \Big\{|m_n(X) -Y|^2 - |m_n(X)- T_{\beta_n} Y|^2 \Big|\mathcal{D}_n \Big \}\\
        && -\textbf{E}\Big\{|m(X) - Y|^2 - |m_{\beta_n}(X)- T_{\beta_n} Y|^2 \Big\}\\
        &=& \textbf{E}\Big\{ (T_{\beta_n} Y -Y)  (2m_n(X)-Y-T_{\beta_n} Y) \Big|\mathcal D_n \Big \}\\
        && - \textbf{E}\Big\{ \Big( (m(X)-m_{\beta_n}(X)) + (T_{\beta_n} Y -Y) \Big) \Big( m(X) + m_{\beta_n}(X) -Y -T_{\beta_n} Y\Big) \Big\}\\
        &=& T_{5,n} + T_{6,n}.
\eeq
With the Cauchy-Schwarz inequality and
\beqm \label{expdurchexp}
I_{\{|Y|>\beta_n\}} \leq \frac{\exp (c_{4}/2 \cdot |Y|^2)}{\exp(c_{4}/2 \cdot \beta_n^2)}
\eeqm
 we conclude
\beq
|T_{5,n}| &\leq& \sqrt{\textbf{E}\big\{|T_{\beta_n} Y -Y|^2 \big\}} \cdot \sqrt{\textbf{E}\big\{|2m_n(X) -Y -T_{\beta_n} Y|^2\big|\mathcal D_n\big\}}\\
        &\leq & \sqrt{ \textbf{E}\big\{ |Y|^2 \cdot I_{\{|Y|>\beta_n\}} \big\}} \cdot \sqrt{\textbf{E}\big\{2\cdot |2m_n(X)
                -T_{\beta_n} Y|^2 + 2\cdot|Y|^2\big|\mathcal D_n\big\}}\\
        &\leq& \sqrt{ \textbf{E}\Bigg\{ |Y|^2 \cdot \frac{\exp(c_{4}/2\cdot|Y|^2)}{\exp(c_{4}/2\cdot \beta_n^2)} \Bigg\}} \\
        && \qquad \cdot
                 \sqrt{\textbf{E}\big\{2\cdot |2m_n(X) -T_{\beta_n} Y|^2 \big|\mathcal D_n \big\}+ 2\textbf{E}\big\{|Y|^2\big\}}\\
                 &\leq& \sqrt{\textbf{E}\Big\{ |Y|^2 \cdot \exp(c_{4}/2\cdot|Y|^2)\Big\}} \cdot
                \exp\left(-\frac{c_{4}\cdot \beta_n^2}{4}\right)
                  \cdot \sqrt{2 (3\beta_n)^2 + 2\textbf{E}\big\{|Y|^2\big\}}.
\eeq
With $ x \leq \exp(x)$ for $x \in \mathbb R$ we get
\[ |Y|^2 \leq \frac{2}{c_{4}} \cdot \exp \left( \frac{c_{4}}{2}\cdot |Y|^2\right) \]
and hence $\textbf{E}\Big\{ |Y|^2 \cdot \exp(c_{4}/2\cdot|Y|^2)\Big\}$ is bounded by
\[
\textbf{E} \left(  \frac{2}{c_{4}} \cdot \exp \left( c_{4}/2 \cdot |Y|^2\right)\cdot \exp(c_{4}/2 \cdot |Y|^2)\right) \leq
\textbf{E} \left(  \frac{2}{c_{4}} \cdot \exp \left( c_{4} \cdot |Y|^2 \right) \right) \leq c_{38}
\]
which is less than infinity by the assumptions of the lemma. Furthermore the third term is bounded by $\sqrt{18 \beta_n^2 + c_{39}}$  because
\beqm \label{Y^2}
\mathbf E(|Y|^2)\leq \mathbf E(1/c_{4} \cdot \exp( c_{4} \cdot |Y|^2) \leq c_{39} < \infty,
\eeqm
 which follows again as above. With the setting $\beta_n = c_{6} \cdot \log(n)  $ it follows for some constants $c_{40}, c_{41} >0$ that
\beq
        |T_{5,n}| &\leq& \sqrt{c_{38}} \cdot \exp \left( -c_{40} \cdot
          \log(n)^2\right) \cdot\sqrt{ (18\cdot c_{6}^2 \cdot (\log n)^2 +c_{39})}
	\leq c_{41} \cdot  \frac{\log(n)}{n}.
\eeq
By the Cauchy-Schwarz inequality we get
\beq
T_{6,n} &\leq& \sqrt{2\cdot \textbf{E}\Bigg\{ | (m(X)-m_{\beta_n}(X))|^2 \Big\}+ 2\cdot  \textbf{E} \Big\{ | (T_{\beta_n} Y -Y) |^2 \Bigg\}}\\
                                && \quad
                         \cdot \sqrt{\textbf {E} \Bigg\{ \Big| m(X) + m_{\beta_n}(X) -Y -T_{\beta_n} Y\Big|^2 \Bigg\}},
\eeq
where we can bound the second factor on the right-hand side in the above inequality in the same way we have bounded the second factor in $T_{5,n}$, because by assumption $||m||_\infty$ is bounded and furthermore $m_{\beta_n}$ is bounded by $\beta_n$. Thus we get for some constant $c_{42}>0$
 \[
 \sqrt{\textbf {E} \Bigg\{ \Big| m(X) + m_{\beta_n}(X) -Y -T_{\beta_n} Y\Big|^2 \Bigg\}} \leq c_{42} \cdot \log(n).
 \]
 Next we consider the first term. By Jensen's inequality it follows that
 \beq
 \textbf{E} \Big\{|m(X)-m_{\beta_n}(X)|^2 \Big\}
 &\leq& \textbf{E} \left\{ \textbf{E} \Big( |Y-T_{\beta_n} Y|^2 \Big| X \Big) \right\} = \textbf{E} \Big\{ |Y-T_{\beta_n} Y|^2 \Big\}.
 \eeq
 Hence we get
 \beq
 T_{6,n} &\leq& \sqrt{ 4 \cdot \textbf{E} \left\{ |Y-T_{\beta_n} Y|^2 \right\}} \cdot c_{42} \cdot \log(n)
 \eeq
 and therefore with the calculations from $T_{5,n}$ it follows that
 $T_{6,n}\leq c_{43} \cdot \log(n)/n$
 for some constant $c_{43}>0$. Altogether we get
 \[T_{1,n} \leq c_{44} \cdot \frac{\log(n)}{n}\]
 for some constant $c_{44}>0$.\\
Next we consider $T_{2,n}$ and conclude for $t>0$
\begin{align*}
\mathbf P\{T_{2,n}>t\}
        &\leq\mathbf P\left\{ \exists f\in T_{\beta_n,supp(X)} \mathcal{F}_n:
        \textbf{E} \left(\left|\frac{f(X)}{\beta_n}-\frac{T_{\beta_n} Y}{\beta_n}\right|^2\right)
        - \textbf{E} \left(\left|\frac{m_{\beta_n}(X)}{\beta_n}-\frac{T_{\beta_n} Y}{\beta_n} \right|^2\right) \right.\\
        & \qquad  -\frac{1}{n} \sum_{i=1}^n \left( \left| \frac{f(X_i)}{\beta_n} - \frac{T_{\beta_n} Y_i}{\beta_n}\right|^2
        -\left|\frac{m_{\beta_n}(X_i)}{\beta_n} - \frac{T_{\beta_n} Y_i}{\beta_n} \right|^2 \right)\\
        & \qquad > \frac{1}{2} \left( \frac{t}{\beta_n^2}
         +\textbf{E} \left(\left|\frac{f(X)}{\beta_n} -\frac{T_{\beta_n} Y}{\beta_n}\right|^2\right) \left.- \textbf{E} \left(\left| \frac{m_{\beta_n}(X)}{\beta_n}
                -\frac{ T_{\beta_n} Y}{\beta_n} \right|^2\right)\right)\right\},
\end{align*}
where $T_{\beta_n,supp(X)} \mathcal{F}_n$ is defined as $\left\{T_{\beta_n}f\cdot 1_{supp(X)}:f\in \F_n\right\}$.
Theorem 11.4 in Gy\"orfi et al. (2002) and the relation
\[
\mathcal N_1 \left( \delta , \left\{ \frac{1}{\beta_n} g :g \in \mathcal G \right\} , x_1^n \right)
\leq \mathcal N_1 \left( \delta \cdot \beta_n, \mathcal G, x_1^n \right)
\]
for an arbitrary function space $\G$ and $\delta >0$ lead to
\begin{align*}
\PROB \{T_{2,n}>t\}
        &\leq  14\cdot \sup_{x_1^n \in supp(X)^n} \mathcal{N}_1 \left(\frac{t}{80\cdot\beta_n}, \mathcal{F}_n, x_1^n \right)
                \cdot \exp \left( - \frac{n}{5136 \cdot \beta_n^2} \cdot t \right).
\end{align*}
 Since the covering number is decreasing in $t$, we can conclude for $\varepsilon_n\geq \frac{80}{n}$
\begin{align*}
\EXP (T_{2,n}) &\leq \varepsilon_n + \int_{\varepsilon_n}^\infty \PROB \{T_{2,n}>t\} dt\\
&\leq \varepsilon_n + 14\cdot \sup_{x_1^n \in supp(X)^n} \mathcal{N}_1 \left(\frac{1}{n\cdot\beta_n}, \mathcal{F}_n, x_1^n \right)
                \cdot \exp \left( - \frac{n}{5136 \cdot \beta_n^2} \cdot \varepsilon_n \right)\cdot \frac{5136 \cdot \beta_n^2}{n}.
\end{align*}
Choosing
\[
\varepsilon_n=\frac{5136\cdot \beta_n^2}{n}\cdot\log\left(
14
\cdot  \sup_{x_1^n \in supp(X)^n} \mathcal{N}_1 \left(\frac{1}{n\cdot\beta_n}, \mathcal{F}_n, x_1^n \right)\right)
\]
(which satisfies the necessary condition
$\varepsilon_n\geq\frac{80}{n}$ if the constant $c_{6}$ in the
definition of $\beta_n$ is not too small)
minimizes the right-hand side and implies
\begin{align*}
  \EXP (T_{2,n}) &\leq \frac{c_{45}\cdot \log(n)^2\cdot
    \log\left( \sup_{x_1^n \in supp(X)^n}
    \mathcal{N}_1 \left(\frac{1}{n\cdot\beta_n},
\mathcal{F}_n, x_1^n \right)
\right)
  }{n}.
\end{align*}
\noindent By bounding $T_{3,n}$ similarly to $T_{1,n}$ we get
\beq
\mathbf E (T_{3,n}) & \leq & c_{46} \cdot \frac{\log(n)}{n}
\eeq
for some large enough constant $c_{46}>0$ and hence we get in total
\beq
\mathbf E \left( \sum_{i=1}^3 T_{i,n} \right) &\leq&  \frac{c_{47} \cdot \log(n)^2
  \cdot
  \left(
\log\left(\sup_{x_1^n \in supp(X)^n}
\mathcal{N}_1 \left(\frac{1}{n\cdot\beta_n},
\mathcal{F}_n, x_1^n \right)
\right)
+1
\right)
}{n}
\eeq
for some sufficient large constant $c_{47}>0$.

We finish the proof by bounding $T_{4,n}$.
Let $A_n$ be the event, that there exists $i \in \{1,...,n\}$ such that $|Y_i|>\beta_n$ and let $I_{A_n}$ be the indicator function of $A_n$. Then we get
\beq
        \textbf{E} (T_{4,n}) &\leq& 2 \cdot\mathbf E \left( \frac{1}{n} \sum_{i=1}^n |m_n(X_i)-Y_i|^2 \cdot I_{A_n}\right)\\
    && + 2 \cdot \mathbf E  \left( \frac{1}{n} \sum_{i=1}^n |m_n(X_i)-Y_i|^2 \cdot  I_{A_n^c}
        - \frac{1}{n} \sum_{i=1}^n |m(X_i)-Y_i|^2  \right)\\
        &=& 2 \cdot\mathbf E \left( |m_n(X_1)-Y_1|^2 \cdot I_{A_n}\right)\\
    && + 2 \cdot \mathbf E  \left( \frac{1}{n} \sum_{i=1}^n |m_n(X_i)-Y_i|^2 \cdot  I_{A_n^c}
        - \frac{1}{n} \sum_{i=1}^n |m(X_i)-Y_i|^2  \right)\\
       &=& T_{7,n} + T_{8,n}.
\eeq
By the Cauchy-Schwarz inequality we get for $T_{7,n}$
\beq
\frac{1}{2} \cdot T_{7,n}
&\leq& \sqrt{ \mathbf E \left( \left( |m_n(X_1) - Y_1|^2\right)^2 \right)} \cdot  \sqrt{ \mathbf P(A_n)}\\
&\leq& \sqrt{ \mathbf E \left( \left( 2 |m_n(X_1)|^2  + 2 |Y_1|^2\right)^2 \right) } \cdot \sqrt{ n \cdot \mathbf P \{|Y_1| >\beta_n \} }\\
&\leq & \sqrt{ \mathbf E \left(  8 |m_n(X_1)|^4  + 8 |Y_1|^4 \right) } \cdot \sqrt{ n \cdot
                    \frac{\mathbf E \left( \exp(c_{4}\cdot |Y_1|^2)\right)}{\exp( c_{4} \cdot \beta_n^2)} },
\eeq
where the last inequality follows as in the proof of
inequality (\ref{expdurchexp}).
Using $x \leq \exp(x)$ for $x \in \mathbb R$ we get
\beq
\mathbf E \left( |Y|^4 \right) &=&\mathbf E \left( |Y|^2 \cdot |Y|^2\right)
                                                        \leq \mathbf E \left( \frac{2}{c_{4}} \cdot \exp\left( \frac{c_{4}}{2} \cdot |Y|^2\right)
                                                         \cdot \frac{2}{c_{4}} \cdot \exp\left( \frac{c_{4}}{2} \cdot |Y|^2\right) \right)\\
                                                         &=& \frac{4}{c_{4}^2} \cdot \mathbf E \left( \exp\left( c_{4} \cdot |Y|^2\right)\right),
\eeq
which is finite by assumption (\ref{subgaus}) of the lemma.
Furthermore $||m_n||_\infty$ is bounded by $\beta_n$ and therefore the first factor is bounded by
\[
c_{48} \cdot \beta_n^2 = c_{49} \cdot (\log n)^2
\]
for some constant $c_{49}>0$. The second factor is bounded by $1/n$, because by the assumptions of the lemma $\mathbf E\left( \exp\left( c_{4}\cdot |Y_1|^2\right) \right)$ is bounded by some constant $c_{50} <\infty$ and hence
 we get
\beq
\sqrt{ n \cdot \frac{\mathbf E \left( \exp(c_{4}\cdot |Y_1|^2)\right)}{\exp( c_{4} \cdot \beta_n^2)} } &\leq&
        \sqrt{n} \cdot \frac{ \sqrt{c_{49}}}{\sqrt{\exp(c_{4} \cdot \beta_n^2)}}
        \leq \frac{\sqrt{n} \cdot
\sqrt{c_{50}}}{\exp((c_{4} \cdot c_{6}^2 \cdot (\log n)^2)/2)}.
\eeq
Since $\exp( - c \cdot \log(n)^2) = O(n^{-2})$ for any $c>0$, we get altogether
\beq
 T_{7,n}  &\leq& c_{51} \cdot \frac{(\log n)^2 \sqrt{n}}{n^2} \leq
 c_{52} \cdot 
\frac{(\log n)^2}{n}.
\eeq
With the definition of $A_n^c$ and $\tilde m_n$ defined as in the assumptions of this lemma we conclude
\beq
T_{8,n} &\leq& 2 \cdot \mathbf E  \left( \frac{1}{n} \sum_{i=1}^n |\tilde{m}_n(X_i)-Y_i|^2 \cdot  I_{A_n^c}
        - \frac{1}{n} \sum_{i=1}^n |m(X_i)-Y_i|^2  \right)\\
&\leq & 2 \cdot \mathbf E  \left( \frac{1}{n} \sum_{i=1}^n |\tilde{m}_n(X_i)-Y_i|^2
        - \frac{1}{n} \sum_{i=1}^n |m(X_i)-Y_i|^2  \right)\\
&\leq & 2 \cdot\mathbf E \left(
\min_{l \in \Theta_n}
\frac{1}{n} \sum_{i=1}^n | g_n (X_i)-Y_i|^2 + pen_n(g_{n,l})
        - \frac{1}{n} \sum_{i=1}^n |m(X_i)-Y_i|^2
        \right)\\
&\leq & 2 \cdot\mathbf E \Bigg(
        \min_{l \in \Theta_n}
        \EXP \Bigg(
        \frac{1}{n} \sum_{i=1}^n | g_{n,l} (X_i)-Y_i|^2 + pen_n(g_{n,l}) \\
        && \hspace*{2cm}
- \frac{1}{n} \sum_{i=1}^n |m(X_i)-Y_i|^2
\bigg| \bb_1^{(1)}, \dots, \bb_r^{(1)}, \dots,
\bb_1^{(I_n)}, \dots, \bb_r^{(I_n)} \Bigg)
        \Bigg)\\
&\leq & 2 \cdot\mathbf E \left(
        \min_{l \in \Theta_n} 
        \int |g_{n,l}(x)-m(x)|^2 \PROB_X (dx) + pen_n(g_{n,l})
\right)
\eeq
because $|T_\beta z - y| \leq |z-y|$ holds for $|y|\leq \beta$.
Hence
\beq
\mathbf E (T_{4,n}) &\leq& c_{53} \cdot \frac{(\log n)^2}{n}+
2 \cdot \mathbf E \left(  \min_{l \in \Theta_n}
        \int |g_{n,l}(x)-m(x)|^2 \PROB_X (dx) + pen_n(g_{n,l})
        \right)
        \eeq
holds.
Thus the proof of Lemma \ref{le6} is complete.\hfill $\Box$\\

\end{appendix}

\end{document}